\documentclass[english,12pt]{smfart}
\usepackage{amssymb}
\usepackage{amsfonts}
\usepackage{amscd}

\swapnumbers

\def\wedgem{\mathop{\wedge}}

\def\ac{{\overline{\rm ac}}}

\def\11{{\mathbf 1}}

\def\Bbase{A_B}

\mathchardef\alphag="7C0B \mathchardef\betag="7C0C
\mathchardef\gammag="7C0D \mathchardef\deltag="7C0E
\mathchardef\varepsilong="7C22 \mathchardef\varphig="7C27
\mathchardef\psig="7C20 \mathchardef\zetag="7C10
\mathchardef\epsilong="7C0F \mathchardef\rhog="7C1A
\mathchardef\taug="7C1C \mathchardef\upsilong="7C1D
\mathchardef\iotag="7C13 \mathchardef\thetag="7C12
\mathchardef\pig="7C19 \mathchardef\sigmag="7C1B
\mathchardef\etag="7C11 \mathchardef\omegag="7C21
\mathchardef\kappag="7C14 \mathchardef\lambdag="7C15
\mathchardef\mug="7C16 \mathchardef\xig="7C18
\mathchardef\chig="7C1F \mathchardef\nug="7C17
\mathchardef\varthetag="7C23 \mathchardef\varpig="7C24
\mathchardef\varrhog="7C25 \mathchardef\varsigmag="7C26
\mathchardef\Omegag="7C0A \mathchardef\Thetag="7C02
\mathchardef\Sigmag="7C06 \mathchardef\Deltag="7C01
\mathchardef\Phig="7C08 \mathchardef\Gammag="7C00
\mathchardef\Psig="7C09 \mathchardef\Lambdag="7C03
\mathchardef\Xig="7C04 \mathchardef\Pig="7C05
\mathchardef\Upsilong="7C07

\newtheorem{theorem}[subsubsection]{Theorem}
\newtheorem{lem}[subsubsection]{Lemma}
\newtheorem{cor}[subsubsection]{Corollary}
\newtheorem{prop}[subsubsection]{Proposition}

\newtheorem{theoremm}[subsection]{Theorem}
\newtheorem{lemm}[subsection]{Lemma}

\newtheorem{propm}[subsection]{Proposition}

\theoremstyle{definition}
\newtheorem{definition}[subsubsection]{Definition}

\newtheorem{def-prop}[subsubsection]{Proposition-Definition}
\newtheorem{def-theorem}[subsubsection]{Theorem-Definition}
\newtheorem{def-lem}[subsubsection]{Lemma-Definition}

\theoremstyle{remark}
\newtheorem{remark}[subsubsection]{Remark}

\theoremstyle{definition}
\newtheorem{definitionm}[subsection]{Definition}

\newtheorem{examplesm}[subsection]{Examples}

\newtheorem{def-propm}[subsection]{Proposition-Definition}
\newtheorem{def-theoremm}[subsection]{Theorem-Definition}
\newtheorem{def-lemm}[subsection]{Lemma-Definition}

\theoremstyle{remark}
\newtheorem{remarkm}[subsection]{Remark}

\theoremstyle{plain}

\numberwithin{equation}{subsection}

\def\boxit#1#2{\setbox1=\hbox{\kern#1{#2}\kern#1}%
\dimen1=\ht1 \advance\dimen1 by #1 \dimen2=\dp1 \advance\dimen2 by
#1
\setbox1=\hbox{\vrule height\dimen1 depth\dimen2\box1\vrule}%
\setbox1=\vbox{\hrule\box1\hrule}%
\advance\dimen1 by .4pt \ht1=\dimen1 \advance\dimen2 by .4pt
\dp1=\dimen2 \box1\relax}

\renewcommand{\theequation}{\thesubsection.\arabic{equation}}

\def\AA{{\mathbf A}}

\def\NN{{\mathbf N}}

\def\QQ{{\mathbf Q}}
\def\RR{{\mathbf R}}

\def\ZZ{{\mathbf Z}}

\def\cF{{\mathcal F}}

\def\cL{{\mathcal L}}
\def\cM{{\mathcal M}}

\def\cP{{\mathcal P}}

\def\cS{{\mathcal S}}
\def\cT{{\mathcal T}}

\DeclareMathOperator*{\sq}{\square}
\def\Jac{{\rm Jac}}

\def\11{{\mathbf 1}}

\newcommand{\comment}[1]{}
\newcommand{\Ko}{{K^\circ}}

\newcommand{\Kt}{{\widetilde K}}

\newcommand{\ord}{\operatorname{ord}}

\renewcommand{\epsilon}{\varepsilon}
\renewcommand{\phi}{\varphi}
\renewcommand{\emptyset}{\varnothing}

\def\ord{{\rm ord}}
\def\ac{{\overline{\rm ac}}}

\def\rv{{ rv}}

\def\Hens{{\rm Hen}}

\def\LHens{L_\Hens{ }}
\def\LHensalt{L_\Hens^{\rm alt}{ }}

\def\LHensast{L_\Hens^{\ast}{ }}

\def\THens{T_\Hens{}}

\begin{document}

\title[$b$-minimality]{$b$-minimality}

\author{Raf Cluckers}
\address{Katholieke Universiteit Leuven,
Departement wiskunde, Celestijnenlaan 200B, B-3001 Leu\-ven,
Bel\-gium. Current address: \'Ecole Normale Sup\'erieure,
D\'epartement de ma\-th\'e\-ma\-ti\-ques et applications, 45 rue
d'Ulm, 75230 Paris Cedex 05, France} \email{cluckers@ens.fr}
\urladdr{www.wis.kuleuven.be/algebra/Raf/}

\author{Fran\c cois Loeser}

\address{{\'E}cole Normale Sup{\'e}rieure,
D{\'e}partement de math{\'e}matiques et applications, 45 rue d'Ulm,
75230 Paris Cedex 05, France (UMR 8553 du CNRS)}
\email{Francois.Loeser@ens.fr} \urladdr{www.dma.ens.fr/$\sim$loeser/}

\begin{abstract}
We introduce a new notion of tame geometry for structures admitting
an abstract notion of balls. The notion is named $b$-minimality and is based on definable families of points and balls.
We develop a dimension theory and
prove a cell decomposition theorem for $b$-minimal structures. We show that $b$-minimality applies
to the theory of Henselian valued fields of characteristic zero, generalizing work by Denef - Pas \cite{Pas} \cite{Pas2}.
Structures which are $o$-minimal, $v$-minimal, or $p$-minimal
and which satisfy some slight extra conditions are also
$b$-minimal, but $b$-minimality leaves more
room for nontrivial expansions. The $b$-minimal setting is
intended  to be a natural framework for the construction of Euler
characteristics and motivic or $p$-adic integrals. The $b$-minimal cell decomposition is a generalization of concepts of P.~J.~Cohen \cite{Cohen}, J.~Denef \cite{Denef2}, and the link between cell decomposition and integration was first made by Denef \cite{Denef}.
\end{abstract}
%\subjclass{03C64,  11S80 (Primary) 03C65, 03C98, 11U09, 12L12, 14A99 (Secondary)}

\maketitle

\section{Introduction}
Originally introduced by Cohen \cite{Cohen} for real and $p$-adic
fields, cell decomposition techniques were developed by Denef
and Pas as a useful device for the study of $p$-adic integrals
\cite{Denef}\cite{Denef3}\cite{Denef2}\cite{Pas}\cite{Pas2}.  Roughly speaking,
the basic idea is to cut a definable set into a finite number of
cells each of which is like a family of balls or points. For general Henselian valued
fields of residue characteristic zero, Denef and Pas proved a cell
decomposition theorem where
the families of balls or points
are parameterized by residue field and value group variables in a definable way.\footnote{In \cite{Cohen}\cite{Denef}\cite{Denef2}\cite{Pas}\cite{Pas2}, cells were very technically defined, and their more simple presentation as families of balls and points appears in \cite{CLoes}.}
An integral over the $p$-adic field is then replaced by a corresponding sum over these residue field and value group variables since the measure of a ball is clear and points have measure zero.
 Denef-Pas cell decomposition plays a fundamental  role in our recent
work \cite{CLoes}  where we lay new general foundations for motivic
integration. When we started in 2002 working on the project that
led to \cite{CLoes}, we originally intended to work in the
framework of an axiomatic cell decomposition of which Denef-Pas cell
decomposition would be an avatar, but we finally decided to keep on
the safe side by staying  within the Denef-Pas framework and we
postponed the axiomatic approach to a later occasion. The present
paper is an attempt to lay the fundamentals of a tame geometry based
upon a cell decomposition into basic families of points and ``abstract'' balls. A key
point in  our approach, already present in  \cite{Pas} and
\cite{CLoes}, is to work in a many sorted language with a unique
main sort and possibly many auxiliary sorts that will parameterize
the families of balls in a definable way. The theory is designed so that no
field structure or topology is required. Instead, only a notion of
balls is needed,  whence the naming $b$-minimality.  The collection
of balls in a model is given by definition as the set of fibers of a
predicate $B$ of the basic language consisting of one symbol $B$.
The notion of $b$-minimality is then based on three axioms, named (b1), (b2) and (b3), where (b1) rather directly imposes the theory to allow cell decomposition and where (b2) and (b3) imply a good dimension theory and exclude pathological behavior.

We show that every $o$-minimal structure is $b$-minimal, but more
exotic expansions of  $o$-minimal structures, like the field of real
numbers with a predicate for the integer powers of $2$ considered by
van den Dries in \cite{vdd2}, are also $b$-minimal, relative to
the right auxiliary sorts. Also $v$-minimal theories of
algebraically closed valued fields, defined by Hrushovski and
Kazhdan in \cite{UdiKazh}, and $p$-minimal theories defined by
Haskell and Macpherson in \cite{Haskell} are $b$-minimal, under some
slight extra conditions for the $p$-minimal case. Our framework is
intended to be versatile enough to encompass promising candidate
expansions, like entire analytic functions on valued and real
fields, but still strong enough to provide cell decomposition and a
nice dimension theory. For $p$-minimality, for example, cell
decomposition is presently missing in  the theory and there seem to be few candidate expansions in sight. For C-minimality and
$v$-minimality, expansions by nontrivial entire analytic functions
are not possible since these have infinitely many zeros in an
algebraically closed valued field. As already indicated, another
goal of the theory is the study of Grothendieck rings and more
specifically the construction of additive Euler characteristics and
motivic
integrals. We intend to go further into that direction in  some future work.\\

Let us briefly review the content of the paper. In section \ref{sbm} the
basic axioms are introduced and discussed and in section \ref{cd}
cell decomposition is proved. Section \ref{sdim} is devoted to
dimension theory. In the next two sections more specific properties
are considered: ``preservation of balls" (which is a consequence of the
Monotonicity Theorem in the $o$-minimal case) and $b$-minimality
with centers (a kind of definable functions approximating the balls in the families). In section \ref{sectionhen} we show that the theory of
Henselian valued fields of characteristic zero is $b$-minimal by
adapting the Cohen - Denef proof scheme. In particular we give (as far
as we know)  the first written instance of cell decomposition in
mixed characteristic for unbounded ramification. Moreover, we  prove
that all definable functions are essentially given by terms. In
section \ref{comp}, we compare $b$-minimality with $p$-minimality,
$v$-minimality and C-minimality. We conclude the paper with some
preliminary results on Grothendieck semirings associated to
$b$-minimal theories.\\

\section{$b$-minimality}\label{sbm}

\subsection{Preliminary conventions} A language may have many
sorts, some of which are called main sorts,
 the others being called auxiliary
sorts. An expansion of a language may introduce new sorts. In this
paper, we shall only consider
languages  admitting a unique main
sort and a distinguished collection of auxiliary sorts. If a model is named $\cM$, then the main sort of $\cM$ is
denoted by $M$.

\par
Usually, for a predicate $B$ of a first order language $\cL$, one
has to specify that it is $n$-ary for some $n$, and one has to fix
sorts $S_1,\ldots,S_n$ such that only for $x_i$ running over $S_i$
one is allowed to write $B(x_1,\ldots,x_n)$. We will not consider
this as information of $\cL$ itself, but instead this information
will be fixed by $\cL$-theories by an axiom of the form
$$
\forall x\  \big(x\in B\ \to  (x\in S_1\times\ldots\times S_n)\big)
$$
for some $n$ and some sorts $S_i$. There is no harm in doing so.

 By \emph{definable} we shall always  mean definable
with parameters, as opposed to $\cL(A)$-definable or $A$-definable,
which means definable with parameters in $A$. By a \emph{point} we
mean a singleton. A definable set is called \emph{auxiliary} if it
is a subset of a finite Cartesian product of (the universes of)
auxiliary sorts.

If $S$ is a sort, then its Cartesian power $S^0$ is considered to be
a point and to be $\emptyset$-definable.
If $S_1$ and $S_2$ are sorts, then by convention $S_1^0$ is the same $\emptyset$-definable point as $S_2^0$, so in our models there is always at least one $\emptyset$-definable point.

Recall that $o$-minimality is about expansions of the language
$\cL_{<}$ with one predicate $<$, with the requirement  that the
predicate $<$ defines a dense linear order without endpoints. In the
present setting we shall  study expansions of a language $\cL_B$
consisting of one predicate $B$, which is nonempty and which has
fibers in the $M$-sort (by definition called balls). In both cases
of tame geometries ($o$-minimal and $b$-minimal), the expansion has to satisfy extra properties.

\subsection{}\label{bmin}
Let $\cL_B$ be the language with one predicate $B$. We require that $B$
is interpreted in any $\cL_B$-model $\cM$ with main sort $M$ as a
nonempty set $B(\cM)$ with
 $$B(\cM)\subset \Bbase \times M$$ where $\Bbase $ is a
finite Cartesian product of (the universes of) some of the sorts of
$\cM$.

When  $a\in \Bbase $ we write $B(a)$ for
$$
B(a):=\{m\in M\mid (a,m)\in B(\cM)\},
$$
and if $B(a)$ is nonempty, we call it a \emph{ball} (in the structure
$\cM$), or $B$-ball when useful.

\begin{definition}\label{bm}
Let $\cL$ be any expansion of $\cL_B$.
We call an $\cL$-model $\cM$ \emph{$b$-minimal} when the following three
conditions are satisfied for every  set of parameters $A$ (the
elements of $A$ can belong to any of the sorts), for every
$A$-definable subset $X$ of $M$, and for every $A$-definable
function $F:X\to M$.
\begin{enumerate}
 \item[(b1)]
There exists an $A$-definable function $ f:X\to S$ with $S$ auxiliary
such that for each $s\in f(X)$ the fiber $f^{-1}(s)$ is a point or a
ball.
 \item[(b2)] If $g$ is a definable function from an auxiliary set to
 a ball, then $g$ is not surjective.
  \item[(b3)] There exists an $A$-definable function $f : X
\rightarrow S$ with $S$ auxiliary such that for each $s\in f(X)$ the
restriction
 $F_{|f^{-1}(s)}$ is either injective or constant.
\end{enumerate}
We call an $\cL$-theory \emph{$b$-minimal} if all its models are
$b$-minimal.
\end{definition}

We call a map $f$ as in (b1) an \emph{$A$-definable
$b$-map on $X$}, or $b$-map for short, and call a map $f$ as in (b3)
\emph{compatible with $F$}.

Typically, in algebraic geometry, one studies families of varieties
by means of a map from one variety $X$ to another variety $Y$, and
then the family consists of the fibers of such a map. Such a family
(of fibers) is then just called a variety over $Y$. Likewise, in
this paper, a definable set $X$ \emph{over} some other definable set
$Y$ is nothing else than a definable map from $X$ to $Y$ and its
interest lies in the study of the family of fibers of this map.
Often, $X$ will be a definable subset of $M^n\times Y$ and the
definable map to $Y$ will be the coordinate projection. This
terminology should not be confused with $X$ being $Y$-definable
meaning that $X$ is definable using parameters from $Y$. By further
analogy, a definable map between two definable sets $X_1$ and $X_2$
over $Y$ is nothing else than a definable function from $X_1$ to
$X_2$ making commutative diagrams with the maps from $X_1$ and $X_2$
to $Y$ which are implicitly understood when we say that the $X_i$
are over $Y$.

With this terminology of definable sets \emph{over} definable sets, we can define a relative version of $b$-maps, that is, $b$-maps over some given definable set.
Namely,
 when $Y$ and $X\subset M\times Y$ are
$A$-definable sets, an \emph{($A$-definable) $b$-map on $X$ over $Y$}
is by definition a ($A$-definable) function $ f:X\to S\times Y $
which makes a commutative diagram with both  projections to $Y$ and
with $S$ auxiliary, such that for each $a\in f(X)$ the fiber
$f^{-1}(a)$ is a point or a ball,  where each $f^{-1}(a)$ is
naturally considered as a subset of $M$
(namely, by projecting on the $M$-coordinate).
 \emph{Compatibility} of $f$
with some definable function $F:X\subset M\times Y\to M\times Y$ over $Y$ is defined
similarly.

\begin{remark}\label{rb2}
Axiom (b2) expresses in a strong way that the auxiliary sorts are different from the main sort.
In particular, the main sort is not interpretable in the auxiliary
sorts by a definable function, in the sense of model theory.
 It follows from (b2) that if $
f:X\to S $ is a $b$-map and $X$ is a ball, then at least one fiber
of  $f$ is not a point, and hence, no ball is a point.
\end{remark}

\begin{remark}
Most of the theory can be
developed with the slightly weaker axiom  (b3$'$)  instead of axiom
(b3):
\begin{itemize}
\item[(b3$'$)]
There exists a finite partition $X_i$ of $X$ into
$A$-definable parts and $A$-definable functions $f_i : X_i
\rightarrow S_i$ with $S_i$ auxiliary such that for each $i$ and
$s\in f_i(X_i)$ the restriction
 $F_{|f_i^{-1}(s)}$ is either injective or constant.
 \end{itemize}
 As soon as there are at least two $\emptyset$-definable auxiliary
points, (b3) and (b3$'$) become equivalent.

Another variant with a similar theory would be to allow for more than one predicate $B$ to define the balls, for example, predicates $B_i\subset A_{B_i}\times M$ and a ball being a fiber $B_i(a)\subset M$ for $a\in A_{B_i}$. We don't pursue this variant here.
\end{remark}

\begin{remark}\label{noncomplete}
Unlike many other notions of minimality or of tame geometries, the notion of $b$-minimality does not require the theory to be complete.
\end{remark}
\subsection{Refinements}
As explained after Definition  \ref{bm}, a definable function $f$ can be seen as the family of its fibers. The following notion of refinement captures the idea that each of the fibers of $f$ gets partitioned into fibers of another function $f'$, where the images of $f$ and of $f'$ are supposed to be auxiliary sets. This last condition is assumed in order to exclude trivial maps like $X\to X$.

Let $f :X\to S$ be a definable function on $X\subset M$, with $S$
auxiliary. By a \emph{refinement of $f$} we mean a pair $(f', g)$
with $f':X\to S'$ and $g : f'(X)\rightarrow S$ definable functions
and $S'$ auxiliary, such that $g\circ f' = f$. Since $g$ is uniquely
determined by $f'$, we shall write $f'$ instead of $(f', g)$ for a
refinement of $f$ and $f \geq f'$ if $f'$ is a refinement of $f$.
This induces a structure of partially ordered set on the set of all
definable functions on $X$ of the form $X\to S$ with $S$ auxiliary.

\par
In the relative setting, let $f :X\to S\times Y$ be a definable
function on $X\subset M\times Y$ over $Y$
(the terminology \emph{over $Y$} means
that $f$ makes a commutative diagram with the projections to $Y$, see the discussion below Definition \ref{bm}), with $S$
auxiliary. By a \emph{refinement of $f$} (over $Y$) we mean a pair
$(f', g)$ with $f':X\to S'\times Y $ a definable function on $X$
over $Y$, $g : f'(X) \rightarrow S\times Y$ a definable function
over $Y$ and $S'$ auxiliary, such that $g\circ f' = f$. Since $g$ is
uniquely determined by $f'$, we call $f'$ a refinement of $f$ and
write $f\geq f'$.

\begin{lem}\label{bref}
Let $\cM$ be a model of a $b$-minimal theory. Let $F:X\to X'$ be a
definable function over $Y$ for some subsets $X$ and $X'$ of
$M\times Y$
(that $F$ is over $Y$ means
that $F$ commutes with the projections to $Y$).
Then any two definable functions $f:X\to S\times Y $ and $f':X\to
S'\times Y $ over $Y$ with $S$ and $S'$ auxiliary admit a common
refinement $f'':X\to S''\times Y$ over $Y$, such that $f''$ is
moreover a $b$-map over $Y$ and compatible with $F$. (In particular,
the opposite category of the category associated to the partially
ordered set of such maps on $X$ over $Y$ is filtering.) If moreover
$F$, $f$, and $f'$ are $A$-definable for some $A$, then $f''$ can be
taken $A$-definable.
\end{lem}
\begin{proof}
For each $y\in Y$ there exists a map as in (b3) on the fiber
$\pi^{-1}(y)$, with $\pi:X\to Y$ the projection, and similarly in
all models of the theory. By compactness, one finds a definable
$f_0:X\to S_0\times Y$ over $Y$, with $S_0$ auxiliary, which is
compatible with $F$.
 Now define $$f_1:X\to S\times S'\times S_0\times Y$$ by $$x\mapsto
p(f(x),f'(x),f_0(x))$$ with $$p: S\times S'\times S_0\times Y^3 \to
S\times S'\times S_0\times Y$$ the projection. For each $s\in
S\times S'\times S_0\times Y$ there exists a $b$-map on
$f_1^{-1}(s)$ by (b1), and this holds in
all models of the theory. By compactness, one finds a $b$-map $f'':X\to
S\times S'\times S_0\times S''\times Y$ on $X$ over $Y$, for some
auxiliary $S''$, whose composition with the projection to $S\times
S'\times S_0\times Y$ equals $f_1$. We obtain this way a map $f''$
as required.
\end{proof}

\subsection{Some Criteria}

The following two criteria are consequences of the lemmas and their
corollary below.

\begin{prop}\label{starst}
Let $\cT$ be an $\cL$-theory with all its models satisfying \textup{(b1)} and \textup{(b2)}. Suppose that there are at least two $\emptyset$-definable
auxiliary points. Then $\cT$ is $b$-minimal if and only if for all
models $\cM$ the following statement \textup{($\ast$)} holds in
$\cM$:
 \begin{enumerate}
 \item[($\ast$)] if $F:X\subset M\to Y$ is a definable
surjection with $Y$ a ball, then not all fibers of $F$ contain
balls.
 \end{enumerate}
\end{prop}

\begin{prop}\label{starstc}
Let $\cT$ be an $\cL$-theory with all its models satisfying \textup{(b1)}.
Suppose that there are at least two $\emptyset$-definable
auxiliary points. Then $\cT$ is $b$-minimal if and only if for all
models $\cM$ the following statement \textup{($\dagger$)} holds in
$\cM$:
\begin{enumerate}
 \item[($\dagger$)] If $F:X\subset S\times M\to Y$ is a definable
surjection with $S$ auxiliary and $Y$ a ball, then there exists
$y\in F(X)$ such that $p( F^{-1}(y))$ does not contain a ball,
with $p$ the projection $X\to M$.
 \end{enumerate}
\end{prop}

\begin{lem}\label{star}
Let $\cT$ be an $\cL$-theory with all its models satisfying \textup{(b2)} and \textup{(b3)}. Then conditions \textup{($\ast$)} and \textup{($\dagger$)}
are satisfied for each model $\cM$ of $\cT$.
\end{lem}

\begin{proof}
We first prove ($\ast$). Let $Y$ be a ball and let $F:X\subset M\to
Y$ be a definable surjection. Consider a map $g: X \rightarrow S$
such that,  for every $s$ in $g(X)$, the restriction
$F_{|g^{-1}(s)}$ is either injective or constant, given by (b3).
  Consider the definable subset $S_1$ of $S$ consisting of all points
$s$ such that $F_{|g^{-1}(s)}$ is injective and set $X_1 :=
g^{-1}(S_1)$, $S_2:=g(X)\setminus S_1$ and $X_2:=g^{-1}(S_2)$.
  Suppose that for some $y\in Y$ the set  $F_{|X_1}^{-1}(y)$ contains
a ball $T$. Then the restriction of $g$ to $T$ is injective, which
contradicts (b2). Hence, no set of the form $F_{|X_1}^{-1}(y)$
contains a ball.
 Now, for every $s$ in $S_2$, $g^{-1}(s)$ is contained in a
(unique) fiber of $F_{|X_2}$. It follows that $F$ induces a
definable surjection $S_2 \rightarrow F(X_2)$. If $Y=F(X_2)$, then
we get a contradiction to (b2). If $Y\not =F(X_2)$ then there
exists a point $y\in Y$ with $F^{-1}(y)=F_{|X_1}{}^{-1}(y)$ which does
not contain a ball as shown above.

\par
We now prove \textup{($\dagger$)}, in a similar way. Let $Y$ be a
ball and let $F:X\subset S\times M\to Y$ be a definable surjection
with $S$ auxiliary. By (b3) and by compactness (as in the proof of
Lemma \ref{bref}) we can consider a map $g: X \rightarrow S\times
S'$ over $S$ such that, for every $t$ in $g(X)$, the restriction
$F_{|g^{-1}(t)}$ is either injective or constant.
 For $s$ in $S$, set  $X_s = X \cap (\{s\} \times M)$.
  Consider the definable subset $S_1$ of $g(X)$ consisting of all points
$t$ such that $F_{|g^{-1}(t)}$ is injective and set $X_1 :=
g^{-1}(S_1)$, $S_2:=g(X)\setminus S_1$ and $X_2:=g^{-1}(S_2)$.
Suppose that for some $y\in Y$ the set $p(F_{|X_1}^{-1}(y))$
contains a ball $T_0$, with $p$ the projection to $M$.
Of course we can identify any subset of $X_s$ with a subset of $M$ by projecting on the $M$-coordinate.
When we do so we have the following helpful claim.
\begin{quote}
\textbf{Claim.}
The set $X_s\cap F_{|X_1}^{-1}(y)$
contains a ball $T$ for some $s$ (after identification with a subset
of $M$
by projecting on the $M$-coordinate).
\end{quote}
 We first prove the claim. Denote by  $C_s$ the set $X_s\cap F_{|X_1}^{-1}(y)$. Suppose by contradiction that, for all $s$, the set $C_s$ contains no ball. Apply (b1) to the sets $C_s$ for all $s$. Since $C_s$ contains no balls, we find $b$-maps $f_s:C_s\to S_s$ for some auxiliary sets $S_s$ such that the nonempty fibers of the $f_s$ are points. But then by compactness we find a single $b$-map $f':p(F_{|X_1}^{-1}(y))\to S'$ for some auxiliary $S'$ such that all nonempty fibers of $f'$ are points. Since we have supposed that $p(F_{|X_1}^{-1}(y))$ contains the ball $T_0$, the restriction of $f'$ to $T_0$ gives a definable bijection between $T_0$ and an auxiliary set. This is a contradiction to (b2) and the claim is proven.
By the claim $g$ is injective on $T$, which again gives a contradiction
to (b2). Hence, no set of the form $p(F_{|X_1}^{-1}(y))$ contains a
ball.
 Now, for every $t$ in $S_2$, $g^{-1}(t)$ is contained in a
(unique) fiber of $F_{|X_2}$. It follows that $F$ induces a
definable surjection $S_2 \rightarrow F(X_2)$. If $Y=F(X_2)$, then
we have a contradiction to (b2). If $Y\not =F(X_2)$ then there
exists a point $y\in Y$ with $F^{-1}(y)=F_{|X_1}{}^{-1}(y)$ which does
not contain a ball as shown above.
\end{proof}

\begin{lem}\label{starplus}
Let $\cT$ be an $\cL$-theory with all its models satisfying \textup{(b1)}, \textup{(b2)} and condition
\textup{($\ast$)}. Suppose that there are at least two
$\emptyset$-definable auxiliary points. Then $\cT$ satisfies \textup{(b3)}.
\end{lem}

\begin{proof}
Let $\cM$ be a model of $\cT$ and let $F:X\to Y$ be $A$-definable
with $X,Y$ subsets of $M$. We may work piecewise on $A$-definable
sets since there are at least two $\emptyset$-definable auxiliary
points. Let $Y_1$ be the definable subset of $Y$ consisting of those
$y\in Y$ such that $F^{-1}(y)$ contains a ball. Let $f_1:Y_1\to S_1$
be a $b$-map on $Y_1$. By ($\ast$), all fibers of $f_1$ are
points. Taking $X_1:=F_1^{-1}(Y_1)$ and $f_1':X_1\to S_1:x\mapsto
f_1(F(x))$, we see that $f_1'$ is compatible with $F_{|X_1}$. Hence
we may suppose that $Y_1$ is empty. Let $\Gamma_F$ be the graph of
$F$.
 Take a $b$-map $f_2:\Gamma_F\to Y\times S$
of $\Gamma_F$ over $Y$ (thus not over $X$). Define
$$
f:X\to S: x\mapsto p_S\circ f_2(x,F(x)),
$$
with $p_S$ the projection $Y\times S\to S$. Then clearly $f$ is
compatible with $F|_{X}$. Indeed, all fibers of $f_2$ are points,
thus for $x_1\not = x_2$ in $X$ either $f(x_1)\not=f(x_2)$, or
$F(x_1)\not = F(x_2)$.

The functions $f_i$ and $f$ can clearly
be taken $A$-definable so (b3) follows.
\end{proof}

\begin{lem}\label{starstb}
Let $\cT$ be any $\cL$-theory with all its models satisfying  condition
\textup{($\dagger$)}. Then all models of $\cT$ satisfy \textup{(b2)}.
\end{lem}
\begin{proof}
Suppose by contradiction to (b2) that $Y$ is a ball and that
$$
g:S\to Y
$$
is definable and surjective, with $S$ auxiliary. Let $T$ be a
ball. Then the map
$$
F:X\to Y:(s,t)\mapsto g(s)
$$
with $X=S\times T$ contradicts ($\dagger$).
\end{proof}

\begin{cor}\label{starstbc}
Let $\cT$ be any $\cL$-theory which satisfies \textup{(b1)} and condition
\textup{($\dagger$)} for all its models $\cM$. Suppose that there are at
least two $\emptyset$-definable auxiliary points.  Then $\cT$
satisfies \textup{(b3)} for all its models.
\end{cor}
\begin{proof}
Follows by Lemmas \ref{starplus} and \ref{starstb} by noticing that
$\cM$ satisfies ($\ast$) whenever it satisfies ($\dagger$).
\end{proof}

\section{Cell decomposition}\label{cd}
Let $\cL$ be any expansion of $\cL_B$, as before, and let $\cM$ be an $\cL$-model.
 Cells are defined by induction on the number of variables.

\subsection{Cells}\label{cel} Let $X\subset M$ be definable and
$f:X\to S$ a definable function with $S$ auxiliary. If all fibers of
$f$ are balls, then we call $(X,f)$ a $(1)$-cell with presentation $f$.
If all fibers of $f$ are points, then we call $(X,f)$ a $(0)$-cell with
presentation $f$. For short, we call such $X$ a cell and $f$ its
presentation.

Let $X\subset M^n$ be definable and let $(j_1,\ldots,j_n)$ be in
$\{0,1\}^n$. Let $p:X\to M^{n-1}$ be a coordinate projection. We call
$X$ a $(j_1,\ldots,j_n)$-cell with presentation
$$
f:X\to S
$$
for some auxiliary $S$, if  for each $\hat
x:=(x_1,\ldots,x_{n-1})\in p(X)$, the set
$p^{-1}(\hat x)\subset \{\hat x\} \times M $,
 identified with a subset of $M$ via the projection $\{\hat x\} \times M \to M$,
is a $(j_n)$-cell with presentation
$$
p^{-1}(\hat x)\to S: x_n\mapsto f(\hat x,x_n)
$$
and $p(X)$ is a $(j_1,\ldots,j_{n-1})$-cell with some presentation
$$
f':p(X)\to S'
$$
satisfying  $f'\circ p=p'\circ f$ for some definable $p':S\to S'$.
\footnote{The condition $f'\circ p=p'\circ f$ for some $p'$ could as
well be left out from the definition of cells since $f$ can always
be refined to imply the existence of $f'$ and $p'$ with this
property, see the proof of \ref{existb}. We chose to include this
condition so that a presentation $f$ captures information about $f'$
(and so on) as well. The same remark applies to the definition of
relative cells and of $b$-maps in \ref{bmap}.}

\subsection{Relative cells} \label{rcel}
In the relative setting, when $Y$ and $X\subset M^n\times Y$ are
definable sets, we say that $X$ together with a definable function $
f:X\to S\times Y$ commuting with the projections $\pi:X\to Y$ and
$S\times Y\to Y$ and with $S$ auxiliary is a $(j_1,\ldots,j_n)$-cell
over $Y$ with presentation
$$
f:X\to S\times Y
$$
if the following holds with $p:M^n\times Y\to M^{n-1}\times Y$ a
coordinate projection. For each $(\hat
x,y):=(x_1,\ldots,x_{n-1},y)\in p(X)$, the set $p^{-1}(\hat
x,y)\subset \{\hat x\} \times M \times \{y\}$,
 identified with a subset of $M$ via the projection $\{\hat x\} \times M \times \{y\} \to M$,
is a $(j_n)$-cell with presentation
$$
p^{-1}(\hat x,y)\to S: x_n\mapsto f(\hat x,x_n,y)
$$
and $p(X)$ is a $(j_1,\ldots,j_{n-1})$-cell over $Y$ with some
presentation
$$
f':p(X)\to S'\times Y
$$
satisfying  $f'\circ p=p'\circ f$ for some $p':S\times Y \to
S'\times Y$.

\subsection{$b$-maps}\label{bmap}
Let $X\subset M^n$ and $f:X\to S$ be definable with $S$ auxiliary.
By induction on the variables, with $p:X\to M^{n-1}$ the coordinate
projection on the first $n-1$ variables, $f$ is called a $b$-map on $X$ when for each $\hat
x:=(x_1,\ldots,x_{n-1})\in p(X)$, the function
 $$p^{-1}(\hat
x)\to S:x_n\mapsto f(\hat x,x_n)
 $$
is a $b$-map on $p^{-1}(\hat x)$ as in section \ref{bmin}, and there
exists some $b$-map
$$
f':p(X)\to S'
$$
satisfying  $f'\circ p=p'\circ f$ for some $p':S\to S'$.

Working relatively, for $X\subset M^n\times Y$ a definable set, we
say that a definable function $ f:X\to S\times Y$ over $Y$ is a
$b$-map on $X$ over $Y$ if there is a projection $p:X\to
M^{n-1}\times Y$ such that for every $(\hat x,y)\in p(X)$ the
restriction of $f$ to $p^{-1}(\hat x,y)$ (also here identified with
a subset of $M$) is a $b$-map on $p^{-1}(\hat x,y)$ and there is a
$b$-map $f':p(X)\to S'\times Y$ on $p(X)$ over $Y$ and a definable
function $p':S\times Y\to S'\times Y$ satisfying $f'\circ p=p'\circ
f$.

\begin{remarkm}
The ordering of coordinates on $M^n$ used for cells and $b$-maps, is
usually implicitly chosen. Such a choice appears also in the
definitions of $o$-minimal and $p$-adic cells.
\end{remarkm}

\begin{def-lemm}[types of cells]\label{type} Let $\cM$ be a model of a $b$-minimal theory.
Let $Y$ and
$X\subset M^n\times Y$ be definable sets.
If $X$ is a $(i_1,\ldots,i_n)$-cell over $Y$, then $X$ is not a
$(i'_1,\ldots,i'_n)$-cell over $Y$ (for the same ordering of the
factors of $M^n$) for any tuple $(i'_1,\ldots,i'_n)$ different from
$(i_1,\ldots,i_n)$. We call $(i_1,\ldots,i_n)$ the type of the cell
$X$.
\end{def-lemm}
\begin{proof}
By induction on $n$. For $n=1$, this follows from (b2), cf. Remark
\ref{rb2}.  The image $X'$  of $X$ under the projection $p_n
:M^n\times Y\to M^{n-1}\times Y$ is a $(i_1,\ldots,i_{n-1})$-cell
over $Y$ and by induction this type is unique. Assume now $X$ is
at the same time a $(i_1,\ldots,i_{n-1},0)$-cell over $Y$ and a
$(i_1,\ldots,i_{n-1},1)$-cell over $Y$. This means that $X$ is at
the same time a $(1)$-cell and a $(0)$-cell over $M^{n-1}\times Y$
which is impossible again by (b2).
\end{proof}

\begin{def-lemm}[Refinements]\label{existb}
Let $\cM$ be a model of a $b$-minimal theory. Let $Y$ and $X\subset
M^n\times Y$ be definable. Then there exists a $b$-map on $X$ over
$Y$. Moreover, any two $b$-maps $f:X\to S\times Y$, $f':X\to
S'\times Y$ over $Y$ have a common refinement, namely, a $b$-map
$f'':X\to S''\times Y$ over $Y$ with (automatically unique)
definable maps $\lambda : f''(X)\to S\times Y$ and $\mu : f''(X)\to
S'\times Y$ such that $\lambda \circ f'' = f$ and $\mu \circ f'' =
f'$.
\end{def-lemm}
\begin{proof}
By compactness (as in the proof of Lemma \ref{bref}) and induction
on $n$ (as in the proof of Lemma \ref{type}).
Indeed, for $n=1$, the existence of a $b$-map on $X$ over $Y$
follows clearly by compactness. For $n>1$, let $f_0:X\to S_0\times
M^{n-1}\times  Y$ be a $b$-map over $M^{n-1}\times  Y$ which exists
by the result for $n=1$. Next, write  $p(X)$ for the image of $X$
under the coordinate projection   $p:M^{n}\times Y\to M^{n-1}\times
Y$. By induction, there exists a $b$-map $f':p(X)\to S'\times Y$
over $Y$. Now let $f:X\to S'\times S_0\times Y$ be the definable
function $x\in X\mapsto (\pi' f'(p(x)),f_0(x))$ with $\pi':S'\times
Y \to S'$ the coordinate projection. Then $f$ is a $b$-map over $Y$
as desired since clearly $f'\circ p=p'\circ f$ with $p':S'\times
S_0\times Y \to S'\times Y$ the coordinate projection. This proves
the existence of $b$-maps on $X$ over $Y$.  The construction of the
refinements is done as in the proof of Lemma \ref{bref}.
\end{proof}

\begin{theoremm}[Cell decomposition]\label{ncd}
Let $\cM$ be a model of a $b$-minimal theory. Let $Y$ and $X\subset
M^n\times Y$ be definable sets. Then there exists a finite partition
of $X$ into cells over $Y$.
\end{theoremm}
\begin{proof}
Same proof as for Lemma \ref{existb}.
\end{proof}

\begin{definitionm}[Refinements of cell decompositions] Let $Y$ and $X\subset M^n\times Y$ be
definable sets. Let $\cP$ and $\cP'$  be two finite partitions of
$X$ into cells $(X_i,f_i)$, resp.~$(X_j',f_j')$, over $Y$. We call
$\cP'$ a \emph{refinement of $\cP$} when for each $j$ there exists
$i$ such that
$$
X_j'\subset X_i
$$
and such that $f_j'$ is a refinement of $f_{i|X_j'}$ in the sense of
Lemma-Definition \ref{existb}, or in other words, for each $b\in
f'_j(X'_j)$, there exists a (necessarily unique) $a\in
f_{i|X_j'}(X'_j)$ such that
$$
  f'_j{}^{-1}(b)\subset f_i^{-1}(a).
$$
\end{definitionm}

\begin{lemm}\label{brefc}
Let $\cM$ be a model of a $b$-minimal theory. Let $Y$ and $X\subset
M^n\times Y$ be definable sets. Then any two cell decompositions of
$X$ over $Y$ admit a common refinement.
\end{lemm}
\begin{proof}
As for Lemma \ref{bref}.
\end{proof}

\begin{remarkm}\label{b1b2}
In fact, the results of this section on cell decomposition already hold for a theory
satisfying only  (b1) and (b2) for all its models.
\end{remarkm}

\section{Dimension theory}\label{sdim}

We now develop a dimension theory for $b$-minimal structures along similar lines as what is done for
$o$-minimal theories, cf.  \cite{vdD}. In what follows $\cL$ is any
expansion of $\cL_B$ as before and $\cM$ is an $\cL$-model.

\begin{definitionm}
 The dimension of a nonempty
definable set $X\subset M^n$ is defined as the maximum of all
sums
$$i_1+\ldots +i_n$$
where $(i_1,\ldots,i_n)$ runs over the types of all cells contained
in $X$, for all orderings of the $n$ factors of $M^n$. To the empty
set we assign the dimension $-\infty$.

If $X\subset S\times M^{n}$ is definable with $S$ auxiliary, the
dimension of $X$ is defined as the dimension of $p(X)$ with
$p:S\times M^{n}\to M^n$ the projection.

When $F:X\to Y$ is an $\cL$-definable function, the dimension of
$X$ over $Y$ is defined as the maximum of the dimensions of the
fibers $F^{-1}(y)$ over all $y\in Y$.
(Of course, the dimension of $X$ over $Y$ depends on $F$.)

We write $\dim (X/Y)$ for the dimension of $X$ over $Y$, and $\dim (X)$
for the dimension of $X$.  The dimension of $X$ over $Y$ is also called the relative dimension of $X$ over $Y$ (along $F$). Usually, $F$ is implicit and $X$ is just called a definable set over $Y$, see the discussion below Definition \ref{bm}.
\end{definitionm}

\begin{propm}\label{ld}
Let $\cM$ be a model of a  $b$-minimal theory. Let $Y,W,Z$ be
definable sets, let $X$ be a $(i_1,\ldots,i_n)$-cell over $Y$, and
let $A,B,C$ be definable sets over $Y$ with $A,B \subset C$. Then
\begin{enumerate}
\item[(0)]  $\dim (X/Y) = i_1+\ldots+i_n$,
\item[(1)]
$ \dim (A\cup B/Y)=\max(\dim(A/Y),\ \dim(B/Y))$,
\item[(2)] $\dim (W\times Z)=\dim(W)+\dim(Z)$.
\end{enumerate}
\end{propm}

\begin{proof}
Let us prove (0). First we notice that for any definable subset
$E\subset X$ which is a $(i_{E_1},\ldots,i_{E_n})$-cell over $Y$,
with the same ordering of coordinates, one has $i_{E_j}\leq i_j$ for
$j=1,\ldots,n$. For $n=1$, this follows from (b2), and for $n>1$
this is proven by induction on $n$ similarly as in the proof of
Lemma-Definition \ref{type}.

Next we show that for any definable subset $E\subset X$ which is a
$(i_{E_1},\ldots,i_{E_n})$-cell over $Y$ with respect to a different
order of the coordinates on $M^n$, one has $i_{E_1}+\ldots
+i_{E_n}\leq i_{1}+\ldots +i_{n}$. This is clear for $n = 1$, so let
us consider the case $n = 2$. We may assume $Y$ is a point. Assume
first that $i_{E_1} = i_{E_2} = 1$. We want to prove that $i_1 = i_2
= 1$.  By (b2), or rather by Remark \ref{rb2}, one finds $i_1=1$.
Denote by  $p_{2X}:M^2\to M$ and $p_{1X}:M\to M^{0}$ the projections
corresponding to the order of coordinates for the cell $X$, and by
$p_{2E}:M^2\to M$ and $p_{1E}:M\to M^{0}$ the projections
corresponding to the order of coordinates for the cell $E$. We may
assume that the image of $E$ by $p_{2E}$ is a ball $T$ and that all
fibers of the restriction of $p_{2E}$
 to $E$
contain
a ball, since $E$
is a $(1,1)$-cell. Assume now that $i_2 = 0$. This means that $X$
is a $0$-cell over $M$ with respect to the projection $p_{2X}$, thus, there
exists an injective map $g : X \rightarrow S\times M$ over $M$,
with $S$ auxiliary. Hence, the map $$F:=p_{2E}\circ
g_{|E}{}^{-1}:g(E)\to M
$$
gives a surjection from $g(E)$ to the ball $T$. Moreover, $p
F^{-1}(t)$ contains a ball for each $t\in T$ with $p:S\times M\to
M$ the projection, which contradicts condition ($\dagger$) of
Lemma \ref{starstc}.

To conclude the case $n = 2$, it is enough to prove that if $i_1 =
i_2= 0$, then $i_{E_1} = i_{E_2} = 0$. But if $i_1 = i_2= 0$, then
$X$ is definably isomorphic to a definable subset of some auxiliary
sorts, which makes it impossible for $E$ to contain a ball by (b2),
hence forces $i_{E_1} = i_{E_2} = 0$. Now for general $n$, it is
enough to consider a transposition
$((x_1,\ldots,x_j,x_{j+1}),\ldots,x_n)\mapsto
((x_1,\ldots,x_{j+1},x_{j}),\ldots,x_n)$ of two adjacent
coordinates. By induction on $n$ and by projecting onto the first
$j+1$ coordinates, $x_1,\ldots,x_{j+1}$, one may suppose that $j+1=n$
and one reduces to the cases already considered. Statement (0)
follows.

Proving (1) amounts to showing that if $X$ is a
$(i_1,\ldots,i_n)$-cell over $Y$ and $X_j$ is a finite partition
of $X$ into $(i_{j1},\ldots,i_{jn})$-cells with respect to the
same ordering of the coordinates, then $\max_j
(i_{j1}+\ldots+i_{jn})=\dim X$, which is clear when $n=1$ and
follows by induction on $n$ when $n>1$. Property (2) is clear by the previous properties since
partitions of $W$ and $Z$ into cells induce a partition  of
$W\times Z$ into cells.
\end{proof}

Recall that, in our terminology, a definable set over another definable set has the meaning as explained after Definition \ref{bm}.

\begin{propm}\label{dimension}
Let $\cM$ be a model of a $b$-minimal theory. Let $Y$  be a
definable set, let $X$ and $X'$ be definable sets over $Y$, and
let $f:X\to X'$ be a definable function over $Y$, that is, compatible with the maps to $Y$. Then
\begin{enumerate}
\item[(3)] $\dim (X)\geq \dim (f(X))$, hence also $\dim (X/Y)\geq \dim
(f(X)/Y)$.
\item[(4)] For each integer $d\geq 0$ the set $S_f(d):=\{x'\in
X'\mid \dim(f^{-1}(x')/Y)=d\}$ is definable and
$$
\dim (f^{-1}(S_f(d))/Y) = \dim (S_f(d)/Y) + d,
$$
with the convention $-\infty + d = -\infty$.
 \item[(5)] If $Y$ is auxiliary, then $\dim (X/Y)=\dim
(X)$.
\end{enumerate}
\end{propm}
\begin{proof}
We first prove (5). We reduce to the case that $X$ is a definable subset of $M^n\times Y$, as follows. By the definition of relative dimensions we may replace $X$ be the graph of $X\to Y$ so that $X$ becomes a definable subset of $M^n\times S\times Y$ for some auxiliary $S$ and some $n\geq 0$. Again by the definition of relative dimensions, we may replace $Y$ by $S\times Y$ to conclude our reduction to the case that $X$ is a definable subset of $M^n \times Y$.
By property (1) of Proposition \ref{ld} and Theorem \ref{ncd} we may then suppose that $X$ is a
$(i_1,\ldots,i_n)$-cell over $Y$.
 Now (5) follows similarly to the way that (0)
of Proposition \ref{ld} is proven. Namely, if $n=1$, (5) follows
from (b2), for $n=2$ it follows from property ($\dagger$), and for
$n>2$ one uses induction.

\par
For (3) and (4), we may suppose that $Y$ is a point, since
relative dimension over $Y$ is defined as the maximum of the
dimensions of the fibers of $y\in Y$.

\par
For (4), let $\Gamma_f\subset X'\times X$ be the graph of $f$
(more precisely, the transpose of the graph).
We have $X\subset M^n\times S$ and $X'\subset M^m\times S'$ for some auxiliary sets $S$ and $S'$ and some $m,n\geq 0$.
We first prove the property when $S$ and $S'$ are singletons, that is, when $X\subset M^n$ and $X'\subset
M^m$.
 By Proposition
\ref{ld} and by Theorem \ref{ncd} we may suppose that $\Gamma_f$ is a
$(i_1,\ldots,i_m,i_{m+1},\ldots,i_{m+n})$-cell.
 Suppose first that $i_{m+1}+\ldots+i_{m+n}=d$. Then $X'=S_f(d)$,
$X=f^{-1}(S_f(d))$, $i_{1}+\ldots+i_{m}=\dim (X')$, and
$i_{1}+\ldots+i_{m+n}=\dim (\Gamma_f)=\dim (X)$ by Proposition
\ref{ld}. Hence (4) follows. When $i_{m+1}+\ldots+i_{m+n}\not=d$
the set $S_f(d)$ is empty and there is nothing to prove.
The case of general $S$ now follows from (5) and compactness. Finally, the case of general $S'$ follows from compactness.

\par
Statement (3) is a corollary of (4) and of (1) of Proposition
\ref{ld}.
\end{proof}

\section{Preservation of balls}\label{spb}

The following property seems especially useful for Henselian valued
fields in the context of motivic integration. Definition \ref{pb} is
used for the change of variables in one variable in motivic
integration in \cite{CLoes}; in an $o$-minimal structure it is a consequence of the
Monotonicity Theorem.

\begin{definitionm}[Preservation of balls]\label{pb}
Let $\cM$ be a $b$-minimal $\cL$-structure. We say that $\cM$
\emph{preserves balls} if for every set of parameters $A$ and
every $A$-definable function
 $$F:X\subset M\to Y\subset M$$
there is an $A$-definable $b$-map
$$
f:X\to S
$$
such that for each $s\in S$ the set
$$
F(f^{-1}(s))
$$
is either a ball or a point.

If moreover $f$ can be chosen in a way that for every refining
$b$-map $f_1:X\to S_1$ the set
 $$
F(f_1^{-1}(s_1)) $$ is also either a ball or a point for each
$s_1\in S_1$, then we say that $\cM$ \emph{preserves all balls}.

We say that a $b$-minimal theory \emph{preserves balls}
(resp.~\emph{preserves all balls}) when all its models do.
\end{definitionm}

\begin{remarkm}\label{pbo}
In an $\cL$-model $\cM$ of a $b$-minimal theory which preserves all
balls, one has the following property which is the $b$-minimal
analogue of the Monotonicity Theorem for $o$-minimal structures, and
which can be proven by taking refinements of $b$-maps. For any
$A$-definable function $F:X\to Y$,
 with $X,Y$ subsets of $M$, there is an $A$-definable $b$-map
$
f:X\to S
$
such that for each $s\in S$ the restriction
 $F_{|f^{-1}(s)}$ is either injective or constant and
$ F(f^{-1}(s)) $ is either a ball or a point and such that similar
properties hold for each refinement of $f$.
\end{remarkm}

\section{$b$-minimality with centers}\label{sbc}

This section is about approximating balls (occurring in definable
families) by definable functions. Indeed, sometimes it is useful to
``center'' balls around some ``center'' points, not necessarily
lying ``in'' the ball, but just lying ``close'' to the ball. In our
treatment we are guided by the situation for characteristic zero
Henselian valued fields. In an $o$-minimal context this notion seems
irrelevant. In Henselian valued fields, centers are heavily used for
the classification of semi-algebraic $p$-adic sets in \cite{C}, for
the definition of motivic integrals in \cite{CLoes} and of
exponential motivic integrals in \cite{CLexp}, for the Fubini
Theorem for motivic integrals in \cite{CLoes}, and for finding
estimates for $p$-adic exponential sums in \cite{Cexp}.

\subsection{}\label{secBn}
Define $\cL_B'$ as
$$\cL_B':=\cL_B\cup\{B_n\}_n,$$
the language $\cL_B$ together with predicates $B_n$ for $n$ running
in some index set, to be interpreted in any model in such a way that
$$
B_n\subset \Bbase \times M
$$
for each index $n$, with $\Bbase $ as defined in section \ref{sbm}.

 Let $\cL$ be any extension of $\cL_B'$ and let $\cM$ be an $\cL$-model.
 For $a\in \Bbase $, write
$$
B_n(a):=\{m\in M\mid (a,m)\in B_n\}.
$$

\begin{definitionm}
A point $x$ in $M$ is called a \emph{$B_n$-center of a ball $Y$}
when there is $a\in \Bbase $ such that $Y=B(a)$ and $x\in B_n(a)$. A
point $x$ in $M$ is called a \emph{$B_n$-center of a point $y$} when
$x=y$.

Let $f:X\subset M\to S$ be a $b$-map. A map $c:f(X)\to M$ is called
a \emph{$B_n$-center of $f$} if $c(s)$ is a $B_n$-center of
$f^{-1}(s)$ for each $s\in f(X)$. When there exists such $n$, we
call $c$ a \emph{center} for $f$.

Centers for $b$-maps $f:X\subset Y\times M\to Y\times S$ over $Y$
are defined similarly.

\end{definitionm}

In the context of Henselian valued fields, the index  $n$ for a
$B_n$-center of a ball describes the distance  between the ball and
the center, see section \ref{shens}, where $n>0$ is an integer.

\begin{definitionm}\label{center}Let $\cM$ be a $b$-minimal $\cL$-model.
Say that $\cM$ is \emph{$b$-minimal with $\{B_n\}_n$-centers}, if
every set of parameters $A$ and every $A$-definable $b$-map
$f:X\subset M\to S$ has an $A$-definable refinement $f'$ with an
$A$-definable $B_n$-center for some $n$.

We say an $\cL$-theory is \emph{$b$-minimal with $\{B_n\}_n$-centers}
(or $b$-minimal with centers for short), if all its models are.
\end{definitionm}

The definition of cells can be adapted naturally to a definition for cells with
$\{B_n\}_n$-centers, as follows.

\begin{definitionm}[Cells with $\{B_n\}_n$-centers] Let
$X\subset M\times Y$ be a cell over $Y$ with presentation $ f:X\to S
\times Y.$ A definable function
$$
c:f(X)\to M
$$
is called a $B_n$-center of the cell $(X,f)$ when it is a
$B_n$-center of the $b$-map $f$ over $Y$. The triple $(X,f,c)$ is
called a cell over $Y$ with $B_n$-center $c$ (and presentation
$f$).

 Let $X\subset M^n\times Y$ be a cell over $Y$ with presentation
$f:X\to S \times Y$. Let $p:X\to M^{n-1}\times Y$ be the coordinate
projection.
 By induction, call
 $$
 (X,f,c_1,\ldots,c_n)
 $$
a cell over $Y$ with $(B_{m_1},\ldots,B_{m_n})$-center
$(c_1,\ldots,c_n)$ if $(X,f,c_n)$ is a cell over $M^{n-1}\times Y$
with $B_{m_n}$-center $c_n$ and $(p(X),f',c_1,\ldots,c_{n-1})$ is a
$(j_1,\ldots,j_{n-1})$-cell with
$(B_{m_1},\ldots,B_{m_{n-1}})$-center $(c_1,\ldots,c_{n-1})$ and
some presentation
$$
f':p(X)\to S'
$$
satisfying  $f'\circ p=p'\circ f$ for some $p':S\to S'$.

For short, call $X$ a \emph{cell with $(B_{m_j})_j$-center} (over
$Y$) when there exist such $f$ and $c_i$, or call $X$ a \emph{cell
with $\{B_n\}_n$-center} (over $Y$) when there exists such tuple
$(B_{m_j})_j$.
\end{definitionm}

\begin{theoremm}[Cell decomposition with centers]\label{ncdc}
Let $T$ be a $b$-minimal theory with $\{B_n\}_n$-centers and let
$\cM$ be a model. Let $Y$ and $X\subset M^n\times Y$ be definable
sets. Then there exists a finite partition of $X$ into cells with
$\{B_n\}_n$-centers (over $Y$).
\end{theoremm}
\begin{proof}
As for Theorem \ref{ncd}.
\end{proof}

\begin{remarkm}\label{brefcc}
The analogue of Lemma \ref{brefc} for cells with centers holds, with
the same proof.
\end{remarkm}

\section{Some examples of $b$-minimal structures
}\label{sectionhen}

\subsection{$o$-minimal structures and non $o$-minimal expansions}\label{somin}
Any $o$-minimal structure $R$ admits  a natural $b$-minimal
expansion by taking as main sort $R$ with the induced structure,
 the two point set $\{0,1\}$ as auxiliary sort and two constant
symbols to denote these auxiliary points. A possible interpretation
for $B$ is easy to find, for example,
\begin{equation*}
\begin{split}
B=\{(x,y,m)\in R^2\times R \mid x<m<y &\mbox{ when }x<y,\\ x<m
&\mbox{ when }x=y,\\\quad \mbox{ and} \ m<y &\mbox{ when } x>y
 \},
 \end{split}
\end{equation*}
so that in the $m$ variable one gets all open intervals as fibers of
$B$ above $R^2$. Property (b3) and preservation of all balls is in
this case a corollary of the Monotonicity Theorem for $o$-minimal
structures.

The notion of $b$-minimality leaves much more room for expansions
than the notion of $o$-minimality: some structures on the real
numbers are not $o$-minimal but are naturally $b$-minimal, for
example, the field of real numbers with a predicate for the integer
powers of $2$ is  $b$-minimal by \cite{vdd2} when adding to the
above language the set of integer powers of $2$ as auxiliary sort
and the natural inclusion of it into $\RR$ as function symbol.

Let us make this more precise. Following van den Dries
\cite{vdd2}, let us consider $(\RR, 2^{\ZZ})$ as an ordered field
with the multiplicative group $2^{\ZZ}$ as a distinguished subset,
so $(\RR, 2^{\ZZ})$ is a structure in the language $L$ of ordered
rings augmented with a $1$-variable predicate $A$. In the paper
\cite{vdd2} van den Dries considers the set of axioms $\Sigma$ for
$L$-structures $(R, A)$ expressing that $R$ is a real closed
field, $A$ a multiplicative subgroup of positive elements, that
$2$ is in $A$, no $x$ in $(1, 2)$ belongs to $A$ and that every $x
>0$ lies in $[y, 2y)$ for some $y$ in $A$. He shows that $\Sigma$
axiomatizes the complete theory of $(\RR, 2^{\ZZ})$. This is done
by adding to the language $L$ new $1$-variable predicates $P_i$,
for $i= 1, \dots, n, \dots$,  and a $1$-variable function symbol
$\lambda$ and adding to $\Sigma$ axioms expressing that $P_n (x)$
holds if $x$ is the  $n$-th power of an element in $A$, that
$\lambda (x) = 0$ for $x \leq 0$ and that, for $x >0$, $\lambda
(x)$ lies in $A$ and $\lambda (x) \leq x < 2 \lambda (x)$, and
proving that  the $L^*$-theory $\Sigma^*$ obtained by adding the
new axioms to $\Sigma$ admits elimination of quantifiers.

Let us consider the two sorted language $\cL$ having as  main sort $R=M$
with the language of ordered rings, auxiliary sort $A'$ with the
language of Presburger groups $(0,1,+,-,\leq,\{\equiv_n\}_n)$, two
function symbols $\lambda' : M \rightarrow A'$ and $a : A'
\rightarrow M$, and the ball predicate $B$ as above.

Let $\Sigma'$ be the $\cL$-theory saying that $R$ is a real closed
field, $A'$ is a $\ZZ$-group, $a:A'\to R^\times$ is a homomorphism
of ordered groups, that $2=a(1)$, that every $x\in R$ with $x
>0$ lies in the interval $[a(\lambda'(x)), 2a(\lambda'(x)))$, and that
$\lambda'(x)$ is $0$ for $x\leq 0$. (It follows from $\Sigma'$ that
$a(0)=1$ and, since $a$ preserves the order, there lies no element
of the form $a(y)$ between $1$ and $2$.)

The following proposition is essentially a corollary of the
quantifier elimination result of \cite{vdd2}.
\begin{prop}\label{2vdD}
Use notation from section \ref{somin}. The theory $\Sigma'$ is $b$-minimal, preserves all balls, and
eliminates quantifiers of the sorts $R$ and $A'$.
\end{prop}
\begin{proof}
Axiom (b2) is clear. Axioms (b1) and (b3), quantifier elimination
and preservation of all balls follow from the quantifier elimination
result of \cite{vdd2}. One can also use criteria \ref{starstc} or
\ref{starst} to prove (b3).
\end{proof}

\subsubsection*{}
It is tempting to hope that the field $\RR$ admits nontrivial
$b$-minimal expansions with entire analytic functions other than the
entire exponential function. Results in this direction have been
obtained by C. Miller in \cite{Miller} and independently by Wilkie
in \cite{Wilkie2}, relatively to an auxiliary sort similar to
$2^\ZZ$ ; Miller and Wilkie add functions like $\sin ( \log x)$ on
the positive real line to the real field. A more careful study is
needed.

\subsection{Henselian valued fields of characteristic
zero}\label{shens}
  In this section we prove that the theory of
Henselian valued fields of characteristic zero is $b$-minimal, in a
natural definitional expansion of the valued field language, by
adapting the Cohen - Denef proof scheme of cell decomposition
\cite{Cohen} \cite{Denef2}. We present a shorter and somewhat
different version of that proof, which, we hope, will enable the
reader to see better the main points in the proof, and to recover
results by Cohen, Denef and Pas \cite{Cohen}, \cite{Denef2},
\cite{Pas}, \cite{Pas2}. For a more detailed proof which is much
closer to the original decision procedure by Cohen, see Pas
\cite{Pas2}. (The more detailed proof in \cite{Pas2} is also word
for word adaptable to the unbounded ramified case.) Note that the
shortcut to prove $p$-adic cell decomposition as presented by Denef
in \cite{Denef} does not yet generalize to any other than the
$p$-adic setting.

As far as we know, this is the first written instance
of cell decomposition in mixed characteristic for unbounded
ramification.
(By bounded ramification we mean that there exists an integer $n>0$ such that $\ord (x^n)>\ord(p) $ for any element $x$ of the maximal ideal, unbounded ramification being the negation of this property.)
 We moreover prove that all definable functions are
essentially given by terms. Preservation of all balls follows.

Let $\Hens$ denote  the collection of all Henselian valued fields of
characteristic zero (hence, mixed characteristic, as well as equal
characteristic zero are allowed).

For $K$ in $\Hens$, write $\Ko$ for the valuation ring, $\Gamma_K$
for the value group, $\ord:K^\times\to \Gamma_K$ for the valuation,
$M_K$ for the maximal ideal of $\Ko$, and $\Kt$ for the residue
field.

For $n>0$ an integer, set $nM_K=\{nm\mid m\in M_K\}$ and
consider
the natural group morphism
$$
rv_{n}:K^\times\to
K^\times/1+nM_K
$$
which we extend to  $rv_n:K\to
(K^\times/1+nM_K)\cup\{0\}$ by sending $0$ to $0$.

For every $n>0$ we write $RV_{n}(K)$ for
$$RV_{n}(K):=(K^\times/1+nM_K)\cup \{0\},$$
and we also write $rv$ for $rv_1$ and $RV$ for $RV_1$.

We use the norm notation $|\cdot|$ for the multiplicative norm associated to the additively written $\ord(\cdot)$; any formula with $|\cdot|$ is an abbreviation of the analogous formula with $\ord$ instead of $|\cdot|$.

We define the  family $B(K)$ of balls by
$$
B(K)=\{(a,b,x)\in K^\times \times K^2 \mid |x-b|<|a|\}.
$$
Hence, a ball is by definition any set of the form $B(a,b)=\{x\in
K\mid |x-b|<|a|\}$ with $a$ nonzero.
 For the centers we consider  the family $\{B_{n}\}_{n}$ over integers
$n>0$, with
$$
B_{n}(K):=\{(a,b,x)\in K^\times \times K^2 \mid |x-b|= |n^{-1}a| \,
\}.
$$
A center in $B_n(K)$ for a ball lies at distance $|n^{-1}|$ from
that ball. Such a center is useful to describe the ball using
definable parameters, as is explained in the following remark.

\begin{remark}[Centers]\label{ballsrv}
For any $n>0$, any nonzero $\xi\in RV_n$, and any $h\in K$, the set
\begin{equation}\label{dball}
X:=\{x\in K\mid rv_n(x-h)=\xi\}
\end{equation}
is an open ball of the form
$$
\{x\in K \mid \ord (x-b) > \ord(n(b-h))\}
$$
for any $b\in X$. Often, none of the points $b$ is definable (over a
certain set of parameters or in a certain definable family of balls)
while $h$ and $\xi$ are definable (uniformly in the family). This is
the advantage of the description (\ref{dball}) of the open ball $X$
and justifies the use of centers. Indeed, with the notation of
sections \ref{bmin} and \ref{secBn}, $X$ is the ball $B(n(b-h),b)$
and $h$ is a center for $X$, that is, $h$ lies in $B_n(n(b-h),b)$.
If some theory has definable centers, then there is a definable
choice for $h$ and thus one can still describe the ball $X$ by
equation (\ref{dball}), even if $b$ is not definable. Given $X$ and
$h$, the element $\xi$ is unique in equation (\ref{dball}). The
converse is also true: if $X$ is a ball $B(a,b)$ and $h$ lies in
$B_n(a,b)$, then $X$ can be written as $\{x\in K\mid
rv_n(x-h)=\xi\}$ for a unique $\xi$.
\end{remark}

We consider the following language  $\LHens$:
it consists of the language of rings $(+,-,\cdot,0,1)$ for the valued field sort which is  the main sort, together with function symbols
$rv_{n}$ for integers $n>0$ from the main sort into the $RV_n$ which
are the auxiliary sorts, and the inclusion language as defined below on
the auxiliary sorts.

We denote by  $\THens$  the theory of all fields in $\Hens$ in the
language $\LHens$.

\subsubsection{The inclusion language on the $RV_n$}\label{induced}
Let $K$ be in $\Hens$. For $a_i\in RV_{n_i}(K)$, $b_j\in K$ and for
$f,g$  polynomials over $\ZZ$ in $n+m$ variables, we let the
expression
$$
f(a_1,\ldots,a_{n},b_1,\ldots,b_m)
$$
correspond to
 the following subset of $K$
$$
\{x\in K\mid (\exists y\in K^{n}) \big( f(y,b)=x \wedgem_i
rv_{n_i}(y_i)=a_i\big)\}.
 $$
Thus, the expression $f(a,b)$ stands for a kind of image of the restriction of $f$ to a specific domain.
By an inclusion
\begin{equation}\label{inc}
f(a,b)\subset g(a,b)
\end{equation}
of such expressions, with $a$ and $b$ tuples as above, we shall mean the inclusions
of the corresponding subsets of $K$.

The inclusion language $\cL_{RV}$ on the sorts $RV_{n}$, $n>0$,
consists of the three symbols $+,\cdot,\subset$, interpreted as the
restriction to (all Cartesian products of) the sorts $RV_n$ of the
relation explained in (\ref{inc}). (There are no terms in this
language, only relations of the form $f(x)\subset g(x)$, for $f$ and
$g$ polynomials formed with $+$ and $\cdot$ in
variables $x$ that run over some of the sorts $RV_n$.)

\subsubsection{An alternative language on the $RV_n$ sorts}
Historically,  in \cite{Basarab}, \cite{KuhlBas}, \cite{Scanlon},
other languages were considered  for  the auxiliary sorts $RV_n$ for
elimination of valued field quantifiers. We define a variant
$\cL_{RV}^{\rm alt}$ of these  languages on the $RV_n$ sorts  that
is most closely related to the variant of Scanlon \cite{Scanlon}
(but without the structure of valued $D$-field in the terminology of
\cite{Scanlon}, or, put otherwise, the language of \cite{Scanlon}
with trivial $D$-structure). We show that our language $\cL_{RV}$ on
the $RV_n$ sorts is definitionally equivalent to that alternative
language $\cL_{RV}^{\rm alt}$ in the sense that the same sets are
definable in both languages. We also show that $\LHens$ eliminates
valued field quantifiers. (The same reasoning holds if one would use
the original Basarab language \cite{Basarab}.)

The language $\cL_{RV}^{\rm alt}$ puts on the $RV_n$ the structure of partially ordered multiplicative semi-groups (it is a multiplicative group with an annihilating zero-element, hence a semi-group). The partial order is the one induced by $\ord$ on the valued field. That is, $\ord(a)<\ord(b)$ for $a,b\in RV_n$ if and only if for $a',b'$ in $K$ with $rv_n(a')=a$, $rv_n(b')=b$ one has $\ord(a')<\ord(b')$. Apart from this structure, there are for each $m$ dividing $n$ the natural projection map $RV_n\to RV_m$  (also denoted by $rv_m$ which is harmless since $rv_m$ on $K$ factorizes through $RV_n$ via $rv_n$) and a partial binary function $+_{n,m}$ from $RV_n^2$ into $RV_m$ sending $(a,b)$ to $c$ if and only $c$ is the unique element satisfying $rv_m(c')=c$ for any $a',b',c'$ in $K$ with $rv_n(a')=a$, $rv_n(b')=b$ and $a'+b'=c'$.

Clearly all symbols of $\cL_{RV}^{\rm alt}$ are valued field quantifier free definable in our language $\cL_{RV}$.

The alternative language $\LHensalt$ is then the language of rings
$(+,-,\cdot,0,1)$ for the valued field sort together with the
function symbols $rv_{n}$ for integers $n>0$ and together with the
language $\cL_{RV}^{\rm alt}$. We know by \cite{Scanlon} (here just
with trivial $D$-structure) that $\LHensalt$ eliminates valued field
quantifiers. Hence, also our language $\LHens$ eliminates valued
field quantifiers. Moreover, both $\LHensalt$ and $\LHens$ are
definitional expansions of the language $\cL_0$ which has the ring
language for the valued field sort, the $RV_n$ sorts, and the
function symbols $rv_n$ for $n>0$. We claim that this implies that
also the languages $\cL_{RV}^{\rm alt}$ and $\cL_{RV}$ are
definitionally equivalent. We only have to prove one direction since
the other is already shown. Take any $\cL_{RV}$-definable subset $X$
in the $RV_n$ sorts. This is also definable in the language $\cL_0$,
possibly using valued field quantifiers, because $\LHens$ is a
definitional expansion of $\cL_0$. Hence, $X$ is also definable in
the language $\LHensalt$ because this is a definitional expansion of
$\cL_0$. By elimination of valued field quantifiers in $\LHensalt$
and since the elements $rv_n(k)$ with $k$ integers are definable in
$\cL_{RV}^{\rm alt}$, the set $X$ is $\cL_{RV}^{\rm alt}$-definable.

We have proven the following variant of the results of \cite{Basarab}, \cite{KuhlBas}, \cite{Scanlon}.
\begin{prop}[Elimination of valued field quantifiers]\label{qehens}
The theory $\THens$ admits elimination of valued field quantifiers in the language $\LHens$. To be precise, for any $\LHens$-formula $\varphi$ there exists an $\LHens$-formula $\psi$ without quantifiers running over the valued field so that $\varphi$ and $\psi$ are equivalent over $\THens$, that is, so that $\THens$ proves $\varphi\leftrightarrow \psi$.
\end{prop}

\begin{remark}Recently, Denef \cite{Denefmanu} gave a
new, alternative,  proof of quantifier elimination for Henselian
valued fields based on monomialization
(which is a strong kind of resolution of singularities, see
\cite{Cutk}).
\end{remark}

The main result of this section is the following.

\begin{theorem}\label{cdh}
The theory $\THens$ is $b$-minimal with $\{B_n\}_n$-centers.
Moreover, $\THens$ preserves all balls.
\end{theorem}

Before we come to the proofs of the main theorems \ref{cdh} and \ref{thens} of this section, we establish two technical lemmas which yield a variant of Cohen's proof \cite{Cohen} of cell decomposition.
 The first lemma is a corollary of Hensel's lemma.
\begin{lem}\label{hm} Let $K$ be in
$\Hens$. Let
 $$f(y)=\sum_{i=0}^m a_iy^i
 $$
be a polynomial in $y$ with coefficients in $K$, let $n>0$ an
integer, and let $x_0\not=0$ be in $RV_{n}(K)$.
 Assume that there exist $i_0>0$ and $x\in K$ satisfying the following conditions
 \textup{(\ref{h0})},  \textup{(\ref{h0b})},  \textup{(\ref{h1})}, and  \textup{(\ref{h2})},
\begin{equation}\label{h0}
rv_n(x)=x_0,
\end{equation}
\begin{equation}\label{h0b}
\mbox{$\ \ord (a_{i_0}x^{i_0})$ is minimal among the $\ord
(a_{i}x^{i})$}
\end{equation}
meaning that $
\min_{0\leq i\leq m} \ord (a_{i}x^{i}) = \ord (a_{i_0}x^{i_0})$,
\begin{equation}\label{h1}
\ord (f(x)) > \ord (n^2  a_{i_0}x^{i_0}),
\end{equation}
and
\begin{equation}\label{h2}
\ord (f'(x)) \leq  \ord (n a_{i_0}x^{i_0-1}).
\end{equation}
Then there exists a unique $y_0\in K$ with
\begin{equation}\label{h3}
f(y_0)=0 \ \mbox{ and }\  rv_{n}(y_0)=x_0.
\end{equation}
Furthermore, if one writes
$f(y) = \sum_i b_i (y-y_0)^i$,
then for any $w$ satisfying
\begin{equation}\label{h0bb}
rv_n(w)=x_0,
\end{equation}
one has moreover
\begin{equation}\label{h5}
\ord f(w)= \ord\, b_1 (w-y_0).
\end{equation}
\end{lem}

\begin{proof}
First consider the case where all coefficients $a_i$ lie in the valuation ring $K^\circ$ and where
$a_{i_0}$ and $x$ are units in $K^\circ$. Then conditions (\ref{h1}) and (\ref{h2}) read $\ord(f(x))>\ord(n^2)$ and $\ord(f'(x))\leq \ord(n)$. Thus, the existence of $y_0$ satisfying $y_0\equiv x\bmod n M_K$ and $f(y_0)=0$  follows from Hensel's
lemma. But $y_0\equiv x\bmod n M_K$ is equivalent to $rv_{n}(y_0)=x_0$ since $x$ is a unit and $rv_n(x)=x_0$ and (\ref{h3}) follows.  To prove  (\ref{h5}), let $w$ satisfy (\ref{h0bb}) and write $f(y) = \sum_i b_i (y-y_0)^i$. Write
$$
f(w)=  b_1 (w-y_0) + b_2(w-y_0)^2 + \ldots,
$$
and
$$
f'(w)=  b_1 + 2 b_2(w-y_0) + \ldots,
$$
where the dots represent higher order terms. By (\ref{h0}),  (\ref{h3}) and (\ref{h0bb}) and since $\ord(f'(x))\leq \ord(n)$, it follows that $\ord(f'(w))\leq \ord(n)$ and $\ord(y_0-w) >n$. Hence, from the expression for $f'(w)$, one gets
 $\ord (b_1)\leq  \ord (n)$. Clearly $\ord(b_j)\geq 0$  for each $j$ since the $a_i$ and $y_0$ lie in $K^\circ$.
Now  (\ref{h5}) is clear from the expression for $f(w)$.

The general case follows after changing coordinates. See Pas
\cite{Pas}, Lemma 3.5,
or \cite{Pas2} for explicit change of variables. In short, one replaces the variable $y$ by $z:=y/x_1$ where $x_1$ is an arbitrary but fixed element satisfying (\ref{h0}), (\ref{h0b}), (\ref{h1}) and (\ref{h2}) for $x=x_1$. Now $g(z):=f(y)/a_{i_0}x_1^{i_0}$ is as in the first case of the proof with the candidate zero of $g$ equal to $1$. Hence there exists a unique root $z_0$ of $g(z)$ with $rv_n(z_0)=rv_n(1)$. Moreover, if one writes $g(z) = \sum_i d_i (z-z_0)^i$, then $\ord g(v)= \ord d_1(v-z_0)$ for any $v$ with $rv_n(v)=rv_n(1)$. Thus $y_0:=z_0x_1$ is a root of $f(y)$ with $rv_n(y_0)=rv_n(x_1)=x_0$. If $y_0$ with $f(y_0)=0$ and $rv_n(y_0)=x_0$ were not unique then this would  contradict the uniqueness of $z_0$ and thus (\ref{h3}) follows.
To prove (\ref{h5}) we may suppose that $\ord  a_{i_0}x_1^{i_0}=0$. For $w$ satisfying (\ref{h0bb}) put $v=w/x_1$. Then $rv_n(v)=rv_n(1)$ and we compute
$$
\ord f(w)= \ord   g(v)= \ord  d_1 (v-z_0) = \ord \frac{d_1}{x_1}(w-y_0)  =  \ord b_1(w-y_0)
$$
since $b_1=d_1/x_1$. Indeed,
$$f(y)=\sum_i b_i (y-y_0)^i=  \sum_i x_1^ib_i (y/x_1-y_0/x_1)^i = \sum_i x_1^ib_i (z-z_0)^i =  \sum_i d_i (z-z_0)^i.$$
 This proves (\ref{h5}).
\end{proof}

\begin{definition}[Henselian functions]\label{hm2}
Let $K$ be in $\Hens$. For all integers $m\geq 0$, $n>0$, define the
function
 $$h_{m,n}:K^{m+1}\times RV_{n}(K) \to K$$ as the function sending
 the tuple
$(a_0,\ldots,a_m,x_0)$ with nonzero $x_0$ to $b$ if there exist
$i_0$ and $x$ that satisfy the conditions (\ref{h0}), (\ref{h0b})
(\ref{h1}), and (\ref{h2}) of Lemma \ref{hm} and where $b$ is the
unique element satisfying (\ref{h3}), and sending
$(a_0,\ldots,a_m,x_0)$ to $0$ in all other cases.

Define $\LHensast$ as the union of the language $\LHens$ together
with all the functions $h_{m,n}$.

\end{definition}

A second main result of this section is the following generalization of some results of \cite{CLoes},
\cite{CLR}, \cite{CLip}.
\begin{theorem}[Term structure]\label{thens}
Let $K$ be in $\Hens$. Let $f:X\to K$ be an $\LHens(A)$-definable
function for some set of parameters $A$. Then there exists an
$\LHens(A)$-definable function $g:X\to S$ with $S$ auxiliary such
that
\begin{equation}\label{et}
f(x)=t(x,g(x))
\end{equation}
for each $x\in X$ and where $t$ is an $\LHensast(A)$-term.
\end{theorem}

\begin{remark}
Note that neither $\LHensast$ nor $\LHens$ have a symbol for the
field inverse on $K^\times$. Indeed, the field inverse is not needed
for Theorem \ref{thens} since the term $h_{1,1}(-1,x,\xi)$ with
$rv(x)\xi=1$ yields a field inverse on $K^\times$.
\end{remark}

The following lemma is a one-parameter variant of Theorem 3.1 of \cite{Pas}.
\begin{lem}\label{lcd}
Let $K$ be in $\Hens$ and let $f(y)$ be a
polynomial over $K$ in the $K$-variable $y$ of degree $d$. Write $A$
for the subset of $K$ consisting of the coefficients of $f$.
 Then:
\begin{enumerate}
\item[(i)] There exist an integer $k>0$  and an $\LHens(A)$-definable
$b$-map
$$
\lambda:K\longrightarrow S
$$
with $B_k$-center
$$
c:\lambda(K)\longrightarrow K
$$
such that, if one writes, for $s=\lambda(y)$,
$$
f(y)= \sum_i a_i(s) (y-c(s))^i,
$$
then
$$
\ord f(y) \leq   \min_{i=0}^d  \ord k a_i(s) (y-c(s))^i.
$$
\item[(ii)] One can ensure that
$c$ is given by an $\LHensast
(A)$-term.

\end{enumerate}

\end{lem}
\begin{proof}
We work by induction on $d$. For $d=0$ the statements are trivial,
so suppose $d>0$.
 Let $f'(y)$ be the derivative of $f$ with respect to $y$. Apply the
induction hypothesis to $f'$. This way, we find a $b$-map
$$\lambda_0:K\longrightarrow S_0$$
 with
$B_{k_0}$-center
$$c_0:S_0\longrightarrow K$$
for some $k_0$ which satisfy (i) and (ii) for $f'$.

Since $c_0$ is  an $\LHensast (A)$-term, there are $\LHensast$-terms
$a_i(s)$
 such that for all $y\in K$ and
$s=\lambda_0(y)$ one has
$$
 f(y) = \sum_{i=0}^d a_i(s) (y-c_0(s))^i.$$

Let $y'$ be the function
$$y':K\to K:y\mapsto y-c_0\circ \lambda_0(y)
$$
and write $a_i'$ for the function
$$
a_i':K\to K:y\mapsto a_i(\lambda_0(y)).
$$

\begin{quote}
\textbf{Claim 1.} \textit{The set $\lambda_0^{-1}(s)$ for $s\in
\lambda_0(K)$ is equal to
$$
\{y\in K\mid rv_{k_0}(y'(y))=s_1\},
$$
for $s_1$ in $RV_{k_0}$ only depending on $s$. Hence, $s_1$ depends in an $\LHens (A)$-definable way on $s$.}\end{quote}

\begin{quote}
\textbf{Claim 2.} \textit{For every nonzero multiple $\ell$ of $k_0$, we may
assume that $c_0$ is a $B_\ell$-center (satisfying (i) and (ii) for
$f'$).
}
\end{quote}
Let us prove the claims. Claim 1 follows from the fact that $c_0$
is a $k_0$-center of $\lambda_0$ and the description of centers and balls in Remark \ref{ballsrv}.
Indeed, for each $s$ in $\lambda_0(K)$, either $\lambda_0^{-1}(s)$ is a ball $B(a,b)$ for some $a$, $b$ (in the notation of section \ref{bmin}) and then $c_0(s)$ lies in $B_{k_0}(a,b)$ and thus one can use the description of Remark \ref{ballsrv}, or, $\lambda_0^{-1}(s)$ is a singleton and then $rv_{k_0}(y'(y))=0$.  That $s_1$ depends in an $\LHens (A)$-definable way on $s$ follows since $s_1$ is uniquely determined by $\lambda_0^{-1}(s)$ and $c_0$ as explained in Remark \ref{ballsrv}.

For Claim 2 and $\ell$ a nonzero multiple of $k_0$, replace
 \begin{eqnarray*}
 S_0 & \mbox{ by }&
 RV_\ell(K)\times S_0,\\
 \lambda_0&\mbox{ by }&
   (rv_\ell\circ y',\lambda_0),\\
 c_0 &\mbox{ by }& c_0\circ p,\ \mbox{ and}\\
 k_0 &\mbox{ by }& \ell,
\end{eqnarray*}
with $p$ the projection  $ RV_\ell(K) \times S_0 \to S_0$. Then (i) and (ii) still hold for
$f'$ and the new $\lambda_0$, $c_0$ and $k_0$. Indeed, since the new $\lambda_0$ has as new component function the function $rv_\ell\circ y'$, that $c_0$ is now a $B_\ell$-center follows from a similar description of balls and centers as explained in Remark \ref{ballsrv} and in the proof of Claim 1. This proves the claims.

We may work piecewise on $\LHens(A)$-definable
pieces since terms on each piece can be combined to one term on the
union, cf.~the proof of Theorem \ref{thens} below. We will work on pieces that we denote by $X$.

By Claim 1, the set $\lambda_0^{-1}(s)$ for $s\in\lambda_0(K)$ is a
ball if and only if $rv_{k_0}(y'(y))\not =0$ for $y$ with
$\lambda_0(y)=s$.
 Hence, we may suppose that each piece $X$ in the finite partition is
a cell with presentation $\lambda_{0|X}$ and center
$c_{0|\lambda_{0}(X)}$.
 Moreover, we may focus on a single piece $X$ which is a $(1)$-cell because the case of
a $(0)$-cell is trivial.

By partitioning further, we may suppose that
$$
X\subset \left\{y\in K \mid  \wedgem_{i,j=0}^d\left( \ord\, a'_i(y)
y'(y)^i
 \sq\nolimits_{ij}\ord\, a'_j(y)y'(y)^j \tilde \sq_j +\infty\right)\right\},
$$
for some $(\sq_{ij})\in\{<,>,=\}^{(d+1)^2}$, $(\tilde \sq_j)_j\in
\{<,=\}^{d+1}$, with the convention that $\ord (0) = +\infty$.
The case where $(\sq_{0j})_j $ is $(=,<,\ldots,<)$ is trivial, with $k$ equal to $k_0$.
Hence, we may suppose that there exists $i_0>0$ such
that
 \begin{equation}\label{i0}
+\infty\not= \ord\, \left(a'_{i_0}(y)y'(y)^{i_0} \right)\leq \ord\,
\left(a'_i(y)y'(y)^i\right) \mbox{ for all $y\in X$ and all $i$.}
\end{equation}

By the previous application of the induction hypothesis to $f'$, we have for $y$ in $X$
\begin{equation}\label{k0i0}
\ord f'(y) \leq  \min_{i}\ord (k_0 i  a'_{i}(y)y'(y)^{i}) \leq  \ord (k_0 i_0  a'_{i_0}(y)y'(y)^{i_0}).
\end{equation}

By Claim 2 we may
\begin{equation}\label{divi}
\mbox{replace $k_0$ by }k_0i_0.
\end{equation}
Then, by (\ref{k0i0}) and (\ref{divi}), one has for all $y\in X$ that
\begin{equation}\label{zx2}
\ord (f'(y))\leq \ord (k_0  a'_{i_0}(y)y'(y)^{i_0}).
\end{equation}

 Define the set $B_1\subset X$ by
\begin{equation*}%{gather*}
B_1:=\left\{y\in X  \mid  \forall\, z\in K \ \left(rv_{k_0}(z) =
rv_{k_0}(y'(y))
\, \to\,   \ord f(z) \leq \ord k_0^2  a'_{i_0}(y)z^{i_0}  \right)\right\},\\
 \end{equation*}%{gather*}
and set  $B_{2}:=
X\setminus B_1$. Note that in particular for $y\in B_1$ one has
\begin{equation}\label{zxbis}
\ord (f(y))\leq \ord (k_0^2  a'_{i_0}(y)y'(y)^{i_0}).
\end{equation}

 We may suppose that either $X=B_1$ or
$X=B_2$, since $B_1$ and $B_2$ are clearly
$\LHens(A)$-definable.

\medskip
\noindent \textbf{Case 1. $X=B_{1}.$}\\
By Claim 2 we may replace $k_0$ by $k_0^2$. Now (i) for $f$  follows from (\ref{zxbis}) with $\lambda=\lambda_0$, $c=c_0$, and $k=k_0$ and (ii) is clear by construction.

\medskip
\noindent \textbf{Case 2. $X=B_{2}.$}\\
By the definition of $B_2$ one has
\begin{equation}\label{zx3}
 \ord (f(z))> \ord (k_0^2
a'_{i_0}(y)z^{i_0}) \mbox{ for $y\in X$ and some $z$ with $rv_{k_0}(z) =
rv_{k_0}(y'(y))$.}
\end{equation}

We will need a new center. First replace $\lambda_0$ as is done in the proof of Claim 2, but just with $\ell=k_0$ (and replace $c_0$ and $S_0$ accordingly). Let  $p: S_0\to RV_{k_0}$ be the projection on the $rv_{k_0}\circ y'
$ component of $S_0$,  which exists since we have replaced $\lambda$ according to the proof of Claim 2 with $\ell=k_0$.
Define the $\LHensast(A)$-term $d:S_0\to K$ by
$$
d(s):=c_0(s) + h_{d,k_0}\big(a_0(s),\ldots,a_d(s),p(s)\big),
$$
with notation from Definition \ref{hm2}.
 By Lemma \ref{hm} and by (\ref{zx2}) and (\ref{zx3}),
for each $y\in X$ and $s=\lambda_0(y)$, the element $d(s)$ lies in the ball
$\lambda_0^{-1}(s)$ and $f(d(s))=0$.
 Define
$$S:=S_0\times RV_{1},
$$
$$
\lambda: X\to S: y\mapsto
\big(\lambda_0(y),rv_{1}(y -  d(\lambda_0(y))  )\big),
$$
$$
c:= d\circ \pi,
$$
with $\pi:S\to S_0$ the
projection,
and define $k=1$. Clearly $c$ satisfies (ii). Also, for $y\in X$ and $s_0=\lambda_0(y)$, one has
\begin{equation}\label{yd}
rv_n(y-c_0(s_0))=rv_n(d(s_0)-c_0(s_0)),
\end{equation}
since $d(s_0)$ lies in the ball $\lambda_0^{-1}(s_0)$ which equals
$$
\{z\mid rv_n(z-c_0(s_0))=rv_n(y-c_0(s_0))\},
$$
see Claim 1.
 Let us check that condition (i) for $f$ holds
on $X$ for the present choice of $\lambda$, $c$, and $k$.
 Writing
$$
f(y)= \sum b_i(s)(y-c(s))^i \mbox{
for $y\in X$ and
 $s=\lambda(y)$},
 $$
one has by (\ref{yd}) and by (\ref{h5}) of Lemma \ref{hm2} that
\begin{equation}\label{Tay}
\ord\, f(y) = \ord\, b_1(s) (y-c(s)),
\end{equation}
which proves the last part of (i).
Since $d(s_0)$ lies in $\lambda_0^{-1}(s_0)$, by the description of balls and centers in Remark \ref{ballsrv}, and by the definition of $\lambda$, it
follows that $d$ is a $B_{1}$-center for $\lambda$. Hence, (i) for
$f$ holds on $X$ for this $\lambda$ and $c$.
\end{proof}

\begin{proof}[Proof of Theorem \ref{cdh}]
The proof of $b$-minimality is based on Lemma \ref{lcd}.
In particular the proof of axiom (b1) for $K$ is derived from Lemma \ref{lcd}. Let
$X\subset K$ be $\LHens(A)$-definable for some set of parameters $A$. By Proposition \ref{qehens}, we may suppose that $X$ is
given by an $\LHens(A)$ formula $\varphi$ without valued field
quantifiers. Let $f_j$ be the polynomials appearing in
$\varphi$. We may suppose that in $\varphi$, the polynomials $f_j$
only appear in the form
\begin{equation}\label{fac}
\rv_{m}(f_j)
\end{equation}
for some $m>0$ since the expression $f_j(x)=0$ is
equivalent to $\rv_{m}(f_j(x))=0$.

Apply Lemma \ref{lcd} to each of the polynomials $f_j$ to find numbers $k_j>0$, $b$-maps $\lambda_j$ and centers $c_j$. As for Claim 1 in the proof of Lemma \ref{lcd}, we may replace each of the $k_j$ by $k:=m\cdot \ell$, with $\ell$ the least common multiple of the $k_j$. It then follows by (i) of Lemma \ref{lcd}  that $\rv_{m}(f_j(x))$ factorizes through $\lambda_j$ for each $j$. In other words, $\rv_{m}(f_j(x))$ does only depend on $\lambda_j(x)$.

Define $\lambda$ as the product map of the $\lambda_j$,
that is,
$$
\lambda: K\to \prod_j S_j: x\mapsto (\lambda_j(x))_j
$$
with $S_j$ the image of $\lambda_j$.

Since a finite intersection of balls is a ball, the map $\lambda$ is a $b$-map. Since clearly the characteristic function of $X$ factorizes through $\lambda$, (b1) follows.

We prove (b2). Examining valued field quantifier free formulas in
one valued field variable as above, one finds that if $S\to K$ is a
definable function from an auxiliary set $S$ into $K$, then its
image is finite
and contained in the zero set of a polynomial.
 Property (b2) thus follows.

For (b3) one uses Proposition \ref{starst} and ($\ast$) is clear by
looking at quantifier free formulas in two valued field variables
which give a definable function. Clearly the graph of such a function must lie in an algebraic set of dimension $\leq 1$. This proves the $b$-minimality.

The preservation of all balls is a special case of Theorem
\ref{cdha}, which is derived in \cite{CLip} from Weierstrass division in rings of analytic functions.
\end{proof}
\begin{proof}[Proof of Theorem \ref{thens}]
 Note that it is
enough to work piecewise. Namely, with the functions $h_{1,1}$ one
can make terms which are characteristic functions of $\rv(1)$,
$\rv(2)$, and so on,  hence one can always paste a finite
number of terms on finitely many disjoint pieces together. An example of a characteristic function of $rv(1)$ is the function $h_{1,1}(1,-1,\cdot):RV\to K:a\mapsto h_{1,1}(1,-1,a)$. If one has two terms $t_j$, $j=1,2$ on subsets $B_j$ of some set $C$, one can replace the $B_j$ and $C$ by $B_j\times \{rv(j)\}$ and $C\times \{rv(1),rv(2)\}$ and construct the single term $\sum_j \chi_{rv(j)} t_j$ with $\chi_{rv(j)}$ the characteristic function of $rv(j)$.

Now let $f$ be an $\LHens(A)$-definable function.
In the theory of $b$-minimality, one derives general cell decomposition from property (b1) by compactness, see the proof of Theorem \ref{ncd}. In the proof of Theorem \ref{cdh}, property (b1) is derived from (i) of Lemma
 \ref{lcd}.
 Similarly to the mentioned application of Lemma \ref{lcd} and the proof of general $b$-minimal cell decomposition by compactness,  it follows that a cell
decomposition theorem holds where all the centers are given by
$\LHensast(A)$-terms. That is, an $\LHens(A)$-definable set can be partitioned into $\LHens(A)$-definable cells whose centers are given by $\LHensast(A)$-terms. Partition the graph of $f$ into such cells
to find the desired piecewise terms.
\end{proof}

\subsection{Analytic structure on Henselian valued fields of characteristic zero}
\label{san}
 The search for an expansion of $\THens$ with a nontrivial entire
analytic function is open and challenging. For the real line, Wilkie and Miller
answered this quest with $\exp$ and other entire functions, see
section \ref{somin}.

Working with analytic functions with bounded domain, the following
is proven in \cite{CLip}, where the analytic functions have as
domains products of the valuation ring and the maximal ideal. The
work in \cite{CLip} generalizes and axiomatizes \cite{Ccell} and \cite{CLR}.

\begin{theorem}[\cite{CLip}]\label{cdha}
Let $\cL_{\rm an}$ be any of the analytic expansions of $\LHens$ as
described in the part of \cite{CLip} on Henselian valued fields and
let $\cT_{\rm an}$ be the corresponding $\cL_{\rm an}$-expansion of
$\THens$.

Then the theory $\cT_{\rm an}$ is $b$-minimal with
$\{B_n\}_n$-centers, preserves all balls, and allows elimination of
valued field quantifiers in the language $\cL_{\rm an}$. The
analogue of Theorem \ref{thens} holds for the language $\cL_{\rm
an}^*$, the union of the language $\cL_{\rm an}$ with all the
functions $h_{m,n}$.
\end{theorem}

\section{Comparison with v-minimality,
p-minimality, and C-minimality}\label{comp}

The comparison with $o$-minimality and a generalization are
already worked out in section \ref{somin}.

\begin{propm}\label{pv1}
Let $T$ be a $p$-minimal theory, as defined in \cite{Haskell}. Let
$T'$ be the theory $T$ with as extra auxiliary sorts the residue
field and the value group and the natural maps into them from the
valued field (that is, the residue modulo the maximal ideal on the valuation ring and zero outside it, resp.~the valuation).
Suppose that $T$ has definable Skolem functions.
Let $B$ and $B_n$ be as in section \ref{shens}. Then
$T'$ is $b$-minimal with $\{B_n\}_n$-centers (with respect to the
auxiliary sorts the residue field and the value group).
\end{propm}

\begin{proof}
The statement  follows from the main result of Mourgues \cite{Mourgues} and
the theory of $p$-minimal fields by Haskell and Macpherson \cite{Haskell}, as follows. If $T$ is any $p$-minimal theory (not necessarily having definable Skolem functions), then $T'$ satisfies
properties (b2) and (b3)  by the results of \cite{Haskell} on
dimensions in $p$-minimal fields. Indeed, (b2) follows from Theorem 3.3 of \cite{Haskell} and (b3) follows from Proposition \ref{starst}, Theorem 6.3 of \cite{Haskell} and properties of $\rm{algdim}$ as defined in \cite{Haskell}.
 By the main result of \cite{Mourgues}, for any $p$-minimal theory $T$ with definable Skolem functions, the theory $T'$ has cell decomposition with centers. From this, property (b1) and the center property \ref{center} for $T'$ follow.
 Thus, if $T$ has
definable Skolem functions, then $T'$ satisfies property (b1) of
Definition \ref{bm} and has the center property \ref{center}, by the main result of \cite{Mourgues}.
\end{proof}

\begin{propm}\label{pv2}
Let $T$ be a $v$-minimal theory of algebraically closed valued
fields, as defined in \cite{UdiKazh} (hence the residue
characteristic is zero and the auxiliary sort is $RV$). Let $B$
and $B_n$ be as in section \ref{shens}. Then $T$ is $b$-minimal with
$B_1$-centers and preserves all balls.
\end{propm}

\begin{proof}
Property (b1) and the center property \ref{center}
follows from Lemma 3.31 of \cite{UdiKazh} by noting
that any definable subset of the valued field and main sort $K$ is a finite Boolean combination
of (open or closed) balls and points. Namely, let $X$ be an
$A$-definable subset of $K$. Write $X$ as a finite Boolean
combination of (open or closed) balls and points. By Lemma 3.9 of \cite{UdiKazh} there exists a definable bijection $h$ between any given definable finite set and an auxiliary set.
Let $h$ be such a bijection between these finitely many points in the Boolean combination and an auxiliary set. Extend $h$ on $K$ by zero to some map $h:K\to S$ with $S$ auxiliary, or extend it on $K$ in some trivial definable way to an auxiliary set $S$. Define the
collection $X'$ of closed balls in $K$ as the union of all the
closed balls in this Boolean combination and for each occurring
open ball in the Boolean combination the minimal closed ball
containing this ball. Then $X'$ is an $A$-definable finite set of
closed balls. By Lemma 3.31 of \cite{UdiKazh} one can take a
finite definable set $Y$ such that in each of the closed balls of
$X'$ there lies exactly one point of $Y$.   Let $Y'$ be the (finite) set of all averages of points in $Y$. Let $g:Y'\to R$ be a definable bijection with $R$ auxiliary, which exists by Lemma 3.9 of \cite{UdiKazh}.
Now let $f$ be the definable map $x\mapsto (rv(x-y),h(z),g(y)) $ into the auxiliary set $RV\times S\times R$ where $y\in Y'$ is the average of the points in $Y$ that lie closest to $x$ and $z$ the average of the finitely many distinguished points in the Boolean combination that lie closest to $x$. Then $f$ is a $b$-map as required for (b1) as follows from the description of balls in Remark \ref{ballsrv}.

 Property (b2) is clear from Lemma 3.41 of \cite{UdiKazh}. Property
(b3) follows by criterium \ref{starst} of section \ref{sbm}, and by the additivity property (a property similar to (4) of Proposition \ref{dimension}) of the $VF$-dimension of \cite{UdiKazh} following from Corollary 3.58 of \cite{UdiKazh}.

Preservation of all balls follows from Proposition 5.1 of
\cite{UdiKazh}. To see this, note that by the first sentence of
the proof of Proposition 5.1 of \cite{UdiKazh}, for a finite
$A$-definable equivalence relation $\sim$ on $K$ there exists  an
$A$-definable $f:K\to S$ for some auxiliary $S$ such that $x\sim
y$ if and only if $f(x)=f(y)$.
\end{proof}

\begin{remarkm}[C-minimality]
Let $\cT$ be a $C$-minimal theory (see \cite{Macpherson} and
\cite{HM}) of algebraically closed valued fields of characteristic
zero and let $\cT'$ be the union of $\cT$ with $\THens$, so, the
auxiliary sorts are the $RV_n$. It is not clear to us whether
$\cT'$ automatically satisfies (b3). Proposition 6.1 of \cite{HM}
seems to give a sufficient condition for (b3) via one of the
criteria of section \ref{sbm}. Theorem 3.11 of \cite{HM} goes in
the direction of preservation of all balls for $\cT'$, but is only
local. Axiom (b1) encompasses somehow the lack of Skolem functions, see also
Lemma 6.6 of \cite{HM}. In positive characteristic
already  (b1) is problematic.
\end{remarkm}

\begin{examplesm}
\item[(1)] In the $p$-adic situation, one can use a technical description of $p$-adic cells based on the definition given by Denef in \cite{Denef}, see \cite{Ccell} for such a description. Write $P_3$ for the nonzero cubes in $\QQ_p$ and then $X = P_3\cap \ZZ_p$ is an example of a $p$-adic cell as in \cite{Ccell}. Suppose for simplicity that $p>3$. To see that $X$ is also a cell in the $b$-minimal setting, let $f:X\to \ZZ^2$ be the definable function given by $x\mapsto (\ord(x), \ac(x))$, where $\ac(x)$ is the first nonzero coefficient of $x$ in the $p$-adic expansion $x=\sum_{i\in\NN} a_ip^i$ with $a_i\in\{0,\ldots,p-1\}$.  (This first nonzero coefficient can be controlled piecewise by using predicates $P_k$ of nonzero $k$th powers for well chosen $k>0$, and their cosets, see \cite{Denef2}.) Clearly $f$ is a $b$-map whose fibers are all balls, hence, $X$ is a $b$-minimal (1)-cell with presentation $f$. For a center for $(X,f)$ one can take the zero function $c:f(X)\to \QQ_p:z\mapsto 0$.

    \item[(2)] In an algebraically closed valued field $K$ with residue field of characteristic zero, it is well known that a definable subset of $K$ in the language of valued fields is a (finite) Boolean combination of balls and points, namely by quantifier elimination. We give two examples that in this case such a Boolean combination is also a finite union of $b$-minimal cells. Let $R$ be the valuation ring of $K$, and $t$ an element of the maximal ideal of $R$. Then $X = R\setminus (t)$ is clearly a simple Boolean combination of balls. Let $f:X\to RV$ be the definable function $x\mapsto rv(x)$. Then again $f$ is a $b$-function whose fibers are all balls, hence $X$ is a $b$-minimal (1)-cell with presentation $f$. For the set $Y=R\setminus \{t\}$, one can use the $b$-map $g:Y\to RV:y\mapsto rv(y-t)$ whose fibers are again all balls, hence also $Y$ is a (1)-cell. For a center for $(X,f)$ one can take the zero function, and for a center for $(Y,g)$ one can take the constant function with value $t$, see Remark \ref{ballsrv}.

\end{examplesm}

A further comparison with tame geometries in model theory and some
further open questions  can be found in \cite{Cb}.

\section{Grothendieck semirings}\label{gsr}
Let $T$
be  a $b$-minimal $\cL$-theory and $\cM$ a model.
In this section, inspired by the
treatments in \cite{UdiKazh} and \cite{CLoes},
we set some lines for a general study of Euler
characteristics on $b$-minimal structures.
\subsection{}
Let $\cS$ be a collection of sets, and $\cF_\cS$ a collection of
functions between some of the sets of $\cS$. Write $A\sim_0 B$ for
$A$ and $B$ in $\cS$ if there exists $f:A\to B$ in $\cF$. Let
$\sim_\cF$ be the equivalence relation on $\cS$ generated by the
relation $\sim_0$. When $A\sim_\cF B$ we say
$A$ and $B$ in $\cS$ are isomorphic.

We define the Grothendieck semigroup $K_0^+(\cS,\cF)$ as the quotient
of the free semigroup on isomorphism classes of $\sim_\cF$ divided
out by the relation
$$
[A]+[B]=[A\cup B]
$$
if $A,B\subset C$ are disjoint and $A,B,C,$ and $A\cup B$ belong to
$\cS$.

\subsection{}
Recall that in this section, $T$ is a $b$-minimal
$\cL$-theory and $\cM$ a model.
For $Y$ a definable set, let $\cS(M/Y)$ be the collection of all
$\cL$-definable subsets of $M^n\times Y$ for $n\geq 0$, and
$\cF(M/Y)$ the collection of $\cL$-definable bijections between such
sets which commute with the projections to $Y$. Likewise, let
$\cS(\cM/Y)$, resp. $\cS(\cM_{\rm Aux}/Y)$, be the collection of all
$\cL$-definable subsets of $X\times Y$ for all definable $X$, resp.
for all auxiliary $X$, and $\cF(\cM/Y)$, resp. $\cF(\cM_{\rm
Aux}/Y)$, the collection of $\cL$-definable bijections between such
sets which commute with the projections to $Y$.

Then write

$$K_0^+(M/Y)\mbox{ for }K_0^+(\cS(M/Y),\cF(M/Y)),$$

$$K_0^+(\cM/Y)\mbox{ for }K_0^+(\cS(\cM/Y),\cF(\cM/Y)),$$
and

$$K_0^+({\rm Aux}/Y)\mbox{ for }K_0^+(\cS(\cM_{\rm Aux}/Y),\cF(\cM_{\rm
Aux}/Y)).$$

These semigroups carry a multiplication induced by Cartesian
product, hence they are endowed with a semiring structure.

\subsection{}
Define $K_0^+({\rm Aux}/Y)[\NN]$ as the graded semigroup
$$
\bigoplus_ {i\in \NN}K_0^+({\rm Aux}/Y),
$$
that is, the direct sum of countably many copies of $K_0^+({\rm Aux}/Y)$ indexed by the nonnegative integers $\NN$.
Write $[X][i]$ for
$(0,\ldots,0,[X],0,\ldots)$  if the
entry $[X]$ occurs at the $i$th position in the tuple. This
semigroup $K_0^+({\rm Aux}/Y)[\NN]$ is a semiring where
multiplication of $[A][i]$ and $[B][j]$ with $A,B\in \cS(\cM_{\rm
Aux}/Y)$ and $i,j\in \NN$ is given by $[A\times B][i+j]$. Note
that this semiring is generated as a semigroup by the elements
$[X][i]$ for $X$ in $\cS(\cM_{\rm Aux}/Y)$ and $i\geq 0$.

\subsection{}
Consider a finite partition of a definable set $X\subset M^n\times
Y$ into $(i_{j1},\ldots,i_{jn})$-cells $X_j$ over $Y$, with
presentations
$$
f_j:X_j\to S_j\times Y
$$
over $Y$, and set  $R_j:=f_j(X_j)$.
 To such an $X$  partitioned into cells with presentations $f_j$ we associate the element
$$\chi(X,f_j)_j:=\sum_j [R_j][\sum_\ell i_{j\ell}]$$
in $K_0^+({\rm Aux}/Y)[\NN]$.

Write $J$ for the ideal of $K_0^+({\rm Aux}/Y)[\NN]$ generated by the
relations
$$
\chi(X,f_j)_j=\chi(X,g_k)_k
$$
for all $X$ in $\cS(M)$ and collections $\{g_k\}$ and $\{f_j\}$ of
$b$-maps of cells corresponding to any two cell decompositions of
$X$ over $Y$.

Then there is a natural map
\begin{equation}\label{chi}
\chi: \cS(M/Y) \to K_0^+({\rm Aux}/Y)[\NN]/J
\end{equation}
sending $X$ to the class of $ \chi(X,f_j)_j$ for some presentations
$\{f_j\}$ of a cell decomposition of $X$.

In this formalism, the following kind of change of variables for $\chi$ becomes almost trivial.

\begin{prop}\label{cv}
Suppose that $\cT$ is $b$-minimal. Then $\chi$ factorizes through
the natural projection
$$
\cS(M/Y)\to K_0^+(M/Y).
$$
\end{prop}
\begin{proof}
Clearly $\chi$ is additive with respect to   disjoint union. Let $X$
and $X'$ be subsets of $ M^n\times Y$ and $ M^\ell\times Y$
respectively and suppose that they are isomorphic over $Y$.
 Let $f:X\to X'$ be a definable bijection over $Y$ and let $\Gamma_f$ be its graph.
 Since a cell decomposition of $X$ over $Y$ induces one of $\Gamma_f$, it is
clear that $\chi(\Gamma_f)=\chi(X)$.
 Similarly, $\chi(\Gamma_f)=\chi(X')$ and thus the proposition is
proved.
\end{proof}

\subsection{}
The study of the map $\chi$, in particular its image and its kernel,
and of the ideal $J$ seems to be  quite fundamental in various settings.
Preservation of balls plays an important  role in the study of $J$. In the
valued field setting, variants with isomorphisms in $\cF(M,|\Jac|)$
or in $\cF(M,rv(\Jac))$, consisting of definable bijections which
are $C^1$ and have constant norm or constant $rv$ of the Jacobian seem even
more important.

In a $v$-minimal setting, such a study has been successfully achieved in
\cite{UdiKazh}. In an $o$-minimal setting,  $\chi$
is easily understood by results in \cite{vdD} on the
Euler characteristic and dimension.

\subsection{}
Let $\cT$ be a theory containing $\THens$ which is $b$-minimal with
centers of level $\{B_{n}\}_{n}$,
for which definable functions are piecewise $C^1$, which preserves
all balls (moreover respecting the Jacobian),
which has the same sorts as $\THens$, and which induces no new
structure on the value group in a strong sense.
  Then we
expect that the whole construction of motivic integration of
\cite{CLoes} goes through. Implementation of notation and
constructions in this generality will be given elsewhere. This
should fill in the gap for mixed characteristic\footnote{However,
one can still imagine that by using some ingenious device there is a
way to avoid the auxiliary sorts $RV_{n}$ for $n>1$ in mixed characteristic.}
motivic integration, and allow to implement the analytic framework
by the results of section \ref{sectionhen}. In particular, the terminology of
definable subassignments and the corresponding constructions make
sense for any such $b$-minimal theory.

\subsection*{Acknowledgments}We express our gratitude to Anand Pillay for insisting
that our set-up should be abstract enough not to refer to a  field structure.
We also
thank Jan Denef and Ehud Hrushovski for stimulating discussions
during the preparation of the paper. We thank the referee for useful comments and the Newton Institute for its hospitality.\\

{\small Our work has been partially supported by the project ANR-06-BLAN-0183.
 During the realization of this project, the first author was
a postdoctoral fellow of the Fund for Scientific Research - Flanders
(Belgium) (F.W.O.) and was partially supported by The European Commission -
Marie Curie European Individual Fellowship with contract number HPMF
CT 2005-007121.}

\bibliographystyle{amsplain}
\bibliography{anbib}

\providecommand{\bysame}{\leavevmode\hbox to3em{\hrulefill}\thinspace}
\providecommand{\MR}{\relax\ifhmode\unskip\space\fi MR }
\providecommand{\MRhref}[2]{%
  \href{http://www.ams.org/mathscinet-getitem?mr=#1}{#2}
}
\providecommand{\href}[2]{#2}
\begin{thebibliography}{10}

\bibitem{Basarab}
{\c{S}}.~Basarab, \emph{Relative elimination of quantifiers for {H}enselian
  valued fields}, Ann. Pure Appl. Logic \textbf{53} (1991), no.~1, 51--74.

\bibitem{KuhlBas}
{\c{S}}.~Basarab and F.-V. Kuhlmann, \emph{An isomorphism theorem for
  {H}enselian algebraic extensions of valued fields}, Man. Math. \textbf{77}
  (1992), no.~2-3, 113 -- 126.

\bibitem{Cb}
R.~Cluckers, \emph{An introduction to b-minimality}, to appear in Conference
  Proceedings of Logic Colloquium 2006, Nijmegen, arXiv:math.LO/0610928.

\bibitem{C}
\bysame, \emph{Classification of semialgebraic $p$-adic sets up to
  semialgebraic bijection}, J. Reine Angew. Math. \textbf{540} (2001),
  105--114.

\bibitem{Ccell}
\bysame, \emph{Analytic $p$-adic cell decomposition and integrals}, Trans.
  Amer. Math. Soc. \textbf{356} (2004), no.~4, 1489 -- 1499,
  arXiv:math.NT/0206161.

\bibitem{Cexp}
\bysame, \emph{Multi-variate {I}gusa theory: Decay rates of $p$-adic
  exponential sums}, Int. Math. Res. Not. IMRN (2004), no.~76, 4093--4108.

\bibitem{CLip}
R.~Cluckers and L.~Lipshitz, \emph{Fields with analytic structure}, preprint
  available at http://www.dma.ens.fr/$\sim$cluckers/.

\bibitem{CLR}
R.~Cluckers, L.~Lipshitz, and Z.~Robinson, \emph{Analytic cell decomposition
  and analytic motivic integration}, Ann. Sci. \'{E}cole Norm. Sup. \textbf{39}
  (2006), no.~4, 535--568, arxiv:math.AG/0503722.

\bibitem{CLoes}
R.~Cluckers and F.~Loeser, \emph{\comment{first}{C}onstructible motivic
  functions and motivic integration}, to appear in Inventiones Mathematicae,
  arxiv:math.AG/0410203.

\bibitem{CLexp}
\bysame, \emph{\comment{second}{C}onstructible exponential functions, motivic
  {F}ourier transform and transfer principle}, to appear in Annals of
  Mathematics, arXiv:math.AG/0512022.

\bibitem{Cohen}
P.~J. Cohen, \emph{Decision procedures for real and $p$-adic fields}, Comm.
  Pure Appl. Math. \textbf{22} (1969), 131--151.

\bibitem{Cutk}
S.~D. Cutkosky, \emph{Resolution of singularities}, Graduate Studies in
  Mathematics, vol.~63, Am. Math. Soc., 2004.

\bibitem{Denef}
J.~Denef, \emph{The rationality of the {P}oincar\'e series associated to the
  $p$-adic points on a variety}, Inventiones Mathematicae \textbf{77} (1984),
  1--23.

\bibitem{Denef3}
\bysame, \emph{On the evaluation of certain $p$-adic integrals}, Th\'eorie des
  nombres, S\'emin. Delange-Pisot-Poitou 1983--84, vol.~59, 1985, pp.~25--47.

\bibitem{Denef2}
\bysame, \emph{$p$-adic semialgebraic sets and cell decomposition}, Journal
  f{\"u}r die reine und angewandte Mathematik \textbf{369} (1986), 154--166.

\bibitem{Denefmanu}
\bysame, \emph{Manuscript on quantifier elimination for {H}enselian valued
  fields},  (2006), Private communication.

\bibitem{vdd2}
{L. van den} Dries, \emph{The field of reals with a predicate for the powers of
  two}, Manuscripta Math. \textbf{54} (1985), no.~1-2, 187--195.

\bibitem{vdD}
\bysame, \emph{Tame topology and o-minimal structures}, Lecture note series,
  vol. 248, Cambridge University Press, 1998.

\bibitem{HM}
D.~Haskell and D.~Macpherson, \emph{Cell decompositions of {${\rm C}$}-minimal
  structures}, Ann. Pure Appl. Logic \textbf{66} (1994), no.~2, 113--162.

\bibitem{Haskell}
\bysame, \emph{A version of o-minimality for the $p$-adics}, J. Symbolic Logic
  \textbf{62} (1997), no.~4, 1075--1092.

\bibitem{UdiKazh}
E.~Hrushovski and D.~Kazhdan, \emph{Integration in valued fields}, Algebraic
  geometry and number theory, Progr. Math., vol. 253, Birkh\"auser Boston,
  Boston, MA, 2006, pp.~261--405.

\bibitem{Macpherson}
D.~Macpherson and C.~Steinhorn, \emph{On variants of o-minimality}, Annals of
  Pure and Applied Logic \textbf{79} (1996), no.~2, 165--209.

\bibitem{Miller}
C.~Miller, \emph{Tameness in expansions of the real field}, Logic Colloquium
  '01 (Urbana, IL), Lect. Notes Log., vol.~20, Assoc. Symbol. Logic, 2005,
  pp.~281--316.

\bibitem{Mourgues}
M.-H. Mourgues, \emph{Corps p-minimaux avec fonctions de {S}kolem
  d{\'e}finissables}, S{\'e}minaire de structures alg{\'e}briques
  ordonn{\'e}es, 1999--2000, pr\'epublication de l'\'equipe de logique
  math\'ematique de Paris 7, pp.~1--8.

\bibitem{Pas}
J.~Pas, \emph{Uniform $p$-adic cell decomposition and local zeta functions},
  Journal f\"ur die reine und angewandte Mathematik \textbf{399} (1989),
  137--172.

\bibitem{Pas2}
\bysame, \emph{Cell decomposition and local zeta functions in a tower of
  unramified extensions of a {$p$}-adic field}, Proc. London Math. Soc. (3)
  \textbf{60} (1990), no.~1, 37--67.

\bibitem{Scanlon}
T.~Scanlon, \emph{Valuation theory and its applications volume {II}}, Fields
  Institute Communications Series, ch.~Quantifier elimination for the relative
  {F}robenius, pp.~323 -- 352, {AMS}, {P}rovidence, 2003, {C}onference
  {P}roceedings of the {I}nternational {C}onference on {V}aluation {T}heory
  ({S}askatoon, 1999).

\bibitem{Wilkie2}
A.~J. Wilkie, \emph{Adding a multiplicative group to a polynomially bounded
  structure},  (October 13, 2006), Lecture at the ENS, Paris.

\end{thebibliography}

\end{document}